\documentclass[12pt]{article}

\usepackage{fullpage,amsfonts,amsmath,amsthm,amssymb,mathrsfs,graphicx}

\usepackage{verbatim,url,epstopdf}

   \usepackage[top=2.54cm, bottom=2.54cm, left=2 cm, right=2 cm]{geometry}
  \newcommand{\lab}[1]{\label{#1}}                

 \usepackage[usenames,dvipsnames]{color}

\newcommand{\remove}[1]{}
\newcommand\eqn[1]{(\ref{#1})}

\newcommand{\be}{\begin{equation}}
\newcommand{\bel}[1]{\begin{equation}\lab{#1}\ }
\newcommand{\ee}{\end{equation}}
\newcommand{\bea}{\begin{eqnarray}}
\newcommand{\eea}{\end{eqnarray}}
\newcommand{\bean}{\begin{eqnarray*}}
\newcommand{\eean}{\end{eqnarray*}}

\newtheorem{thm}{Theorem}
\newtheorem{cor}[thm]{Corollary}

\newtheorem{lemma}[thm]{Lemma}

\def\proof{\noindent{\bf Proof.\ }  }
\def\qed{~~\vrule height8pt width4pt depth0pt}

  \newcommand{\nat}{\ensuremath {\mathbb N} }
  \newcommand{\bit}{\clubsuit}

\def\Heavy{{\mathcal H}}
\def\G{{\mathcal G}}
\def\po{{\bf Po}}
\def\ex{{\bf E}}
\def\pr{{\bf P}}

\def\C{{\cal C}}
\def\M{{\cal M}}
\def\MM{{\bf M}}

\def\MMS{{\MM_{\rm simple}}}

\def\P{{\cal P}}

\newcommand{\ctwo}{28}

\def\eps{\epsilon}
\def\la{\lambda}

\def\ss{\smallskip}
\def\non{\nonumber}
\def\no{\noindent}

\date{}

\title{Enumeration of graphs with a heavy-tailed degree sequence}
\author{
Pu Gao\thanks{Research supported by the NSERC PDF. Current address: School of Mathematical Sciences, Monash University, jane.gao@monash.edu}\\
Department of Computer Science\\
 University of Toronto\\
pu.gao@utoronto.ca
\and
Nicholas Wormald\thanks{Research supported by the Canada Research Chairs program and NSERC while  at the Department of Combinatorics \&  Optimization, University of Waterloo. Supported by   Australian Laureate Fellowships grant FL120100125.}
\\
School of Mathematical Sciences\\
Monash University\\
nick.wormald@monash.edu }
\begin{document}
\maketitle
\begin{abstract}

In this paper, we asymptotically enumerate graphs with a given degree sequence  ${\bf d}=(d_1,\ldots,d_n)$ satisfying   restrictions designed to permit heavy-tailed sequences in the sparse case (i.e.\ where the average degree is rather small).
Our general result requires upper bounds on functions of  $M_k= \sum_{i=1}^n [d_i]_k$ for a few small integers $k\ge 1$. Note that $M_1$ is simply   the total degree of the graphs. As special cases, we asymptotically enumerate graphs with
(i) degree sequences satisfying $M_2=o(M_1^{ 9/8})$; (ii) degree sequences following a power law with parameter $\gamma>5/2$; (iii)  power-law degree sequences that mimic   independent power-law ``degrees"  with parameter $\gamma>1+\sqrt{3}\approx 2.732$;  (iv) degree sequences following a certain ``long-tailed'' power law; (v) certain bi-valued sequences. A previous result on sparse graphs by McKay and the second author applies to a wide range of degree sequences but requires $\Delta =o(M_1^{1/3})$, where $\Delta$ is the maximum degree.  Our new result applies in some cases when  $\Delta$ is only barely $o(M_1^ {3/5})$. Case (i) above generalises a result of Janson  which requires $M_2=O(M_1)$ (and hence  $M_1=O(n)$ and $\Delta=O(n^{1/2})$).   Cases (ii) and (iii) provide  the first asymptotic enumeration results applicable to degree sequences of real-world networks following a power law, for which it has been empirically observed   that   $2<\gamma<3$.

\end{abstract}

\section{Introduction}
\lab{s:intro}

For a positive integer $n$, let ${\bf d}=(d_1,d_2,\ldots,d_n)$ be a non-negative integer vector. How many simple graphs are there with degree sequence ${\bf d}$? We denote this number by $g({\bf d})$. This is a natural question, but there is nevertheless no simple formula known for $g({\bf d})$. However, some  simple  formulae have been obtained for the asymptotic behaviour of $g({\bf d})$ as $n\to\infty$, provided certain restrictions are imposed on the degree sequence $\bf d$. Such formulae have been used in many ways, for instance in proving properties of typical graphs with given degree sequence, or for proving properties of other random graphs by classifying them according to degree sequence. They can also lead to new algorithms for generating these graphs uniformly at random.

 We denote $\max_i d_i$  by $\Delta$ throughout this article.  Bender and Canfield~\cite{BC} gave the first general result on the asymptotics of $g({\bf d})$ for the case that
there is a fixed upper bound  on $\Delta$  for all $n$.   This upper bound was later
relaxed by Bollob\'{a}s~\cite{B} to a slowly growing function of $n$. A much more significant relaxation came when
 McKay~\cite{M2} introduced the switching method to this problem (to be explained later in this article), resulting in
 an asymptotic formula when   $\Delta=o(M_1^{1/4})$.  Here and throughout this paper, $M_j=\sum_{i=1}^n  [d_i]_j$  for any integer
 $j\ge 1$, where  $[x]_j$ denotes $x(x-1)\cdots(x-j+1)$.   By improving the switching operations in~\cite{M2}, McKay and Wormald~\cite{MW2}
 further relaxed the constraint on the maximum degree to $o(M_1^{1/3})$.

These results apply best when none of the degrees $d_i$  deviate greatly from the average degree $M_1/n$. But there are important classes of graphs whose degree sequences have ``heavy tails"\!\!. For instance, the degree sequence of the Internet graph exhibits a power-law behaviour~\cite{FFF}, by which it is meant that the number of vertices with degree $k$ is approximately proportional to $k^{-\gamma}$ for some constant $\gamma>1$. Many real networks (e.g.\ web graphs, collaboration networks and many social networks) have such degree sequences and are consequently called scale-free.

 Motivated by research on scale-free networks,  various models have been proposed  to generate random graphs with power-law degree sequences. Some of the well known ones are:  the preferential attachment model~\cite{AB,BRST} and  variations of it~\cite{ACL,BCDR,CF,M}; random hyperbolic graphs~\cite{BP,GPP,K,PKBV,PR}; versions of random graphs with given expected degrees~\cite{CL} (generalising the classical random graph model by Erd\H{o}s and R\'{e}nyi~\cite{ER}); versions of random multigraphs with given degree sequence~\cite{ACL2}.

A natural way to generate random graphs with a given degree distribution  is to sample the degree sequence with the correct distribution, and then generate a random graph with the specified degrees under the uniform probability distribution.  In Section~\ref{sec:results} we define the pairing model, often called the configuration model, which is commonly used to study random graphs with a specified degree sequence.  

Chung and Lu~\cite{ACL2} suggested using the pairing model to generate random graphs with power-law degrees, ignoring what is the main problem with the model: the resulting graph can have loops or multiple edges, i.e.\ is actually a multigraph. Besides the question of whether a multigraph is realistic in modelling real networks, a major problem is that these multigraphs are not uniformly distributed, even though the simple graphs it produces are uniformly distributed.

Producing non-simple graphs is not always an insurmountable problem  for proving properties of the random simple graphs: if the probability that the multigraph is simple is bounded below by a quantity $B$, then the  probability that a simple graph has  some specified (typically undesired) property is at most $1/B$ times the corresponding probability for the multigraph. If $B$ is not too small, the resulting bound can be useful. 
Bollob{\'a}s~\cite{B} instigated this approach, and much subsequent work  used the pairing model to prove properties of the random simple graphs when there is a fixed upper bound on the degrees. In that case,  we may take $B$ to be a positive constant.

It has been  observed empirically that most real-world scale-free networks have degree sequences following a power law whose parameter satisfies $2<\gamma<3$.  Unfortunately, to date there has been no good estimate of $B$ for this important class of degree sequences.  

It is well known that  computing  the probability that the pairing multigraph is simple is equivalent to enumerating graphs with the given degree sequence (see  Section~\ref{sec:results} for more detail).
When $2<\gamma<3$, the power-law graphs have a linear (in $n$) number of edges, whereas the maximum degree is much greater than $n^{1/3}$. Hence, the asymptotic result in~\cite{MW2} cannot be applied. 
A recent result of  Janson~\cite{J,J2} also deals with the case of a linear number of edges,  i.e.\ $M_1=O(n)$, giving an asymptotic approximation for $g({\bf d})$ in the case that $M_2=O(M_1)$. This applies to some cases   not covered by the result from~\cite{MW2}, such as when one vertex has degree approximately $\sqrt n$ and the others have bounded degree.  However,   a power-law degree sequence with $\gamma\le3$ also fails to obey  $M_2=O(M_1)$.

In this paper, we take the next step required for proving properties of graphs with real-world power-law degree sequences, by solving the asymptotic enumeration problem (equivalently, estimating the probability that the pairing model gives a simple graph in such cases) for  a certain range of  values of the parameter $\gamma$ below 3.
Here, since $M_2\ne O(M_1)$, the expected number of loops and multiple edges in the pairing model increases to some power of $n$, suggesting that the probability of a simple graph is likely to be exponentially small.   Estimating this probability with desired precision consequently becomes much more difficult than say in~\cite{J}, where $M_2 = O(M_1)$ and consequently the probability being estimated is bounded away from 0. 

 For estimating such tiny probabilities, a proven approach is to use switchings. Rather than the original switchings used in~\cite{M2} and~\cite{J}, we use the more sophisticated switchings of~\cite{MW2} with refinements that allow us to   control the large error terms caused by  vertices of large degree. The refinements are  necessary because of the difficulty of getting uniform error terms when the degrees of vertices can vary wildly. To this end, we introduce a method in which the multigraphs in the model are classified according to how many  multiple edges join any given pair of vertices. This is a much more elaborate classification structure than has previously been used with switching arguments. 

Our main result applies to a much more general class of degree sequences than  the power-law case, including some that are even more heavy-tailed and some denser. As applications of our general result, we give several interesting examples. In several of these, the maximum degree can significantly exceed $\sqrt{n}$, resulting in some multiple edges with multiplicity tending to infinity. This creates some of the difficulties of the analysis. Indeed, the maximum degree can   reach close to $n^{2/3}$.

Among the examples are two variations of power-law degree sequences. One is along the lines of the model used by  Chung and Lu~\cite{CL}, in which the  number of vertices of degree $k$ is at most $cnk^{-\gamma}$ for some constant $c>0$, uniformly for all integers $k\ge 1$. This  implies the maximum degree is  $O(n^{1/\gamma})$.   We enumerate these graphs when $\gamma>5/2$.  The second version mimics a degree sequence whose components are independent power-law variables,  conditional on even sum. This gives the distribution of  the  maximum degree a longer tail, reaching well above $\sqrt{n}$. For this version, our result requires $\gamma>1+\sqrt{3}\approx 2.732$.  
 

\section{Main results}
\lab{sec:results}

Recall from Section~\ref{s:intro} the definition of  the ``moment'' $M_k$ and maximum degree $\Delta$. We will give an asymptotic estimate of $g({\bf d})$ when $M_2$ (and perhaps also $M_3$ and $M_4$) does not grow too fast compared with $M_1$  without any additional restriction on $\Delta$. We assume throughout this paper that  every component in ${\bf d}$ is at least 1  since results for graphs with vertices of degree 0 then follow trivially. For a valid degree sequence,  $M_1$ must be even. For brevity,  we do not restate this trivial constraint in the hypotheses of our results.
 We use the Landau notation $o$ and $O$.
All asymptotics in this paper refers to $n\to\infty$.

Random graphs with given degree sequence ${\bf d}$ can be generated by the pairing model.
This is a probability space
 consisting of $n$ distinct bins $v_i$ (representing the $n$ vertices), $1\le i\le n$, each containing $d_i$ points, and all points are uniformly at random paired (i.e.\ the points
 are partitioned uniformly at random  subject to  each part containing exactly two points). We call each element in this probability space a {\em pairing}, and two paired points (points contained in the same part) is called a {\em pair}.
  Let $\Phi$ denote the set of all pairings. Then $|\Phi|$ equals the number of matchings on $M_1$ points, and
\be \lab{phi}
|\Phi| = \frac{ M_1!}{2^{M_1/2}(M_1/2)!} =\sqrt 2(M_1/e)^{M_1/2}\big(1+O( M_1^{-1})\big).
\ee
  For each $\P\in\Phi$, let $G(\P)$ denote the multigraph generated by $\P$ by representing bins as vertices and pairs as edges. Thus, $G(\P)$ has degree sequence ${\bf d}$. It is easy to see that every simple graph with degree sequence ${\bf d}$ corresponds to exactly $\prod_{i=1}^n d_i!$ distinct pairings in $\Phi$. Hence,
by letting $\G^*(n,{\bf d})$ denote the probability space of the random multigraphs generated by the pairing model, it follows immediately that
\begin{equation}
g({\bf d})=\frac{|\Phi|}{ \prod_{i=1}^n d_i!} \pr(\G^*(n,{\bf d})\ \mbox{is simple}) \lab{equiv}
\end{equation}
where $|\Phi|$ is given in~\eqn{phi}.
Thus, enumerating graphs with degree sequence ${\bf d}$ is equivalent to estimating the probability that $\G^*(n,{\bf d})$ is simple.

A major difficulty in estimating $\pr(\G^*(n,{\bf d})\ \mbox{is simple})$ using switching arguments is that the vertices with large degrees easily cause big error terms. In order to keep the errors under control, a simple trick is to restrict the maximum degree $\Delta$, such as assuming $\Delta=o(M_1^{1/3})$  as  in~\cite{M2,MW2}. We are able to impose a less severe restriction on the maximum degree by completely reorganising and refining the switching arguments.

Before presenting our general results on the estimates of $g({\bf d})$, or equivalently, $\pr(\G^*(n,{\bf d})\ \mbox{is simple})$, we give several results that are interesting and are simpler. In the following theorem, we consider any degree sequence such that $M_2$ does not grow too fast compared with $M_1$, with   no additional restriction on $\Delta$. 
\begin{thm}\lab{thm:M2}
Let  ${\bf d}$ have minimum component at least 1 and satisfy  $M_2=o(M_1^{9/8})$. Then with $\la_{i,j}=d_id_j/M_1$ and   $|\Phi|$  given  in~\eqn{phi},
\be
g({\bf d})=\big(1+O(\sqrt{\xi})\big) \frac{|\Phi|}{\prod_{i=1}^n d_i!}\exp\left(-\frac{M_1}{2}+\frac{M_2}{2M_1}-\frac{M_3}{3M_1^2}+\frac{3}{4}+\sum_{i<j}\Big(\log(1+\la_{i,j})\Big)\right),\lab{g(d)}
\ee
where $\xi=M_2^4/M_1^{9/2}+M_2^{3/2}/M_1^2+1/M_1$ and necessarily $\xi=o(1)$.
  \end{thm}

\no {\bf Remark}:   Easy calculations show that if $M_3=o(M_1^{3/2})$ then only the first two terms in the expansion of the logarithm contribute to the main terms,  and they produce 
$\sum_{i<j}  \la_{i,j}=  \frac12  M_1  -\frac12 M_2/ M_1-\frac12 +o(1)$ and
$\sum_{i<j} - (1/2)\la_{i,j}^2 =-\frac14   M_2^2/M_1^2-\frac12 M_2/M_1 -\frac14+o(1)$. Hence the main term from the exponential   is $\exp(  -\frac14   M_2^2/M_1^2-\frac12 M_2/M_1) $, which gives    $\pr(\G^*(n,{\bf d})\ \mbox{is simple})$ in agreement with the previously known formulae such as  from~\cite{BC}.
 Moreover, the theorem provides a main term that agrees with~\cite[Theorem 1.4]{J} for  $M_2=O(M_1)$ (which is assumed throughout~\cite{J}), and comes with a much sharper error estimate.  We verify the agreement of these main terms at the end of this section.

 The main advantage of our results over existing ones is for the case when the degree sequence is far from that of a regular graph.  Four  special cases of our main result,  exemplifying this, are given next.

In the first  example, we consider  degree sequences ${\bf d}$ that follow  a so-called power law with parameter $\gamma>1$, i.e.\ the number $n_i$ of vertices of degree $i$ is approximately
$ ci^{-\gamma} n $ for some constant $c>0$.  We relax  the conditions of~\cite{CL} a little and define ${\bf d}={\bf d}(n)$ to be  a {\em power-law density-bounded sequence} with parameter $\gamma$ if  there exists $C>0$ such that  the number of components in ${\bf d}$ taking value $i$ is  at most  $Ci^{-\gamma} n$  for all $i\ge 1$  (and all $n$).

\begin{thm}\lab{thm:powerLaw} Assume that ${\bf d}$  is a power-law density-bounded sequence with parameter $\gamma> 5/2$. Then putting $M_i^*=M_i+M_1$ for $i=2$ and 3, and with $\Phi$ given in~\eqn{phi},
\bean
g({\bf d}) &=& \frac{|\Phi|}{ \prod_{i=1}^n d_i!}\exp\left(-\frac{M_2}{2M_1}-\frac{M_2^2}{4M_1^2}+\frac{M_3^2}{6M_1^3}+\frac{M_4}{4M_1^2}-\frac{M_4^2}{8M_1^4}-\frac{M_6}{6M_1^3}+O\left(\frac{M_2^*\sqrt{M_3^*}}{M_1^2}\right)\right)\\
&=&  \frac{|\Phi|}{ \prod_{i=1}^n d_i!}\exp\left(-\frac{M_2}{2M_1}-\frac{M_2^2}{4M_1^2}+\frac{M_3^2}{6M_1^3}+O(n^{5/\gamma-2})\right)  \quad \mbox{if $5/2<\gamma<3$}.
\eean
\end{thm}

\no {\bf Remarks}. In the case of a ``strict'' power law,  where ${\bf d}$ is a sequence with $n_i= \Theta(  i^{-\gamma} n)$ for $i\le \Delta = \Theta(n^{1/\gamma}) $, with $5/2<\gamma<3$   constant,  the whole exponential factor is $\exp\big(-\Theta(n^{6/\gamma-2})\big)$.  It is difficult to express the exponential factor as a sharper function of $n$ and $\gamma$ alone, since even $M_2^2/M_1^2$ is not sharply determined in this case.

 For the second example, we investigate another type of power-law sequence based on the distribution of a sequence of $n$ independent  identically distributed (i.i.d.) power-law variables, conditional on even sum. Define $F_{\gamma}(i)=\sum_{j\ge i} j^{-\gamma}$. We say a sequence is {\em power-law distribution-bounded with parameter $\gamma$} if 
there exists $C>0$ such that the number of components in the sequence taking value at least $i$ is  at most $CF_{\gamma}(i)n$ for all $i$ and $n$. These are similar to power-law density-bounded sequences except that the upper tail of the distribution,  where $Cni^{-\gamma}$, the ``expected'' number of components equal to $i$, falls below 1, extends further. This results in the maximum component having order $n^{1/(\gamma-1)}$ (instead of $n^{1/\gamma}$).  Our general theorem yields the following enumeration result for such degree sequences. We focus only on $\gamma<3$ as the case $\gamma>3$ follows easily from previously known results (e.g.~\cite{J}) and the case $\gamma=3$ can be easily worked out by applying our general theorem.
\begin{thm}\lab{thm:iid}
Assume that ${\bf d}$  is a power-law distribution-bounded sequence with parameter $3>\gamma>1+\sqrt{3}\approx 2.732$. Then putting $M_i^*=M_i+M_1$ for $i=2$ and 3, and with $\Phi$ given in~\eqn{phi},
\bean
g({\bf d}) &=& \frac{|\Phi|}{\Delta!^{\ell}\delta!^{(n-\ell)}}\exp\bigg(-\frac{M_1}{2}+\frac{M_2}{2M_1}+\frac{3}{4}+\sum_{i<j} \log(1+d_id_j/M_1) +O\big(\xi\big)\bigg),
\eean
where $\xi=n^{(2+2\gamma-\gamma^2)/(\gamma-1)}$.
\end{thm}

Note that   $\xi = o(1)$ for the values of $\gamma$ under consideration in the theorem. Note also that if the entries of ${\bf d}$ are chosen by i.i.d.\ power-law variables conditioned on even sum, with parameter $\gamma'>1+\sqrt 3$, then with probability tending to 1 as $n\to\infty$,  {\bf d} is power-law distribution-bounded with any parameter $\gamma<\gamma'$. If we additionally ensure $\gamma> 1+\sqrt 3$,   the theorem applies.

Next we consider degree sequences with only two distinct degrees, which we call {\em bi-valued}.

\begin{thm}\lab{thm:bidegrees}
Let $3 \le  \delta\le\Delta$ be integers depending on $n$, and assume that   $d_i\in\{\delta,\Delta\}$ for $1\le i\le n$. Let $\ell$ denote the number of vertices with degree $\Delta$. If
\begin{enumerate}
\item[(a)]$\Delta=O(\sqrt{\delta n+\Delta \ell})$  and  $\xi: =  (\Delta^7\ell^3+\Delta^3\delta^4 n^2\ell +\delta^7n^3)/(\delta^4n^4+\Delta^4\ell^4)=o(1)$, or
 \item[(b)] $\Delta=\Omega( \sqrt{\delta n})$ and
$\displaystyle
\xi:=\frac{\Delta^5\ell^3}{\delta^3n^3}+\frac{\Delta^5\ell^2}{\delta^2n^3}+\frac{\delta^3}{n}+\frac{\Delta^3\ell }{n^2}=o(1)$,
\end{enumerate}
then
\be
g({\bf d})=\frac{|\Phi|}{\Delta!^{\ell}\delta!^{(n-\ell)}}\exp\bigg(-\frac{M_1}{2}+\frac{M_2}{2M_1}+\frac{3}{4}+\sum_{i<j} \log(1+d_id_j/M_1) +O\big(\sqrt{\xi}\big)\bigg)\lab{bivalue}
\ee
where  $\Phi$ is given in~\eqn{phi} and $M_i$ is simply $[\Delta]_i\ell+[\delta]_i(n-\ell)$. \end{thm}

\no {\bf Remarks}.

\no (i)    For convenience we omit the cases $\delta=1$ and 2, which can be worked out easily from our main result but require a   different statement.

\no (ii)
 The summation in the exponent in~\eqn{bivalue} is  easy to express in terms of $\delta$ etc.\ as there are only three possible values of $d_id_j$.

\noindent (iii) If we apply Theorem~\ref{thm:bidegrees}(a) with $\Delta=\delta=d$ and $\ell=n$, then we obtain the asymptotic formula for the number of $d$-regular graphs for $d=o(n^{1/3})$, which agrees with~\cite{M2} (note that the regular case is extreme in the opposite direction from what we are aiming at here, which is highly irregular degree sequences).  

\noindent (iv)
Theorem~\ref{thm:bidegrees}(b)  applies to some instances of bi-valued sequences where the minimum degree is around $ n^{1/3-\eps}$ and simultaneously there are up to $n^{\eps}$ vertices with maximum degree as large as $n^{2/3-\eps}$. These are much higher degrees than can be reached by any previously published results on enumeration of sparse graphs with given degree sequence.

\no (iv) The bi-valued degree sequence  easily generalises to a much wider class of degree sequences as follows. The first $\ell$ vertices have degree at most $\Delta$; for some $\la>0$, $\sum_{\ell<j\le n}d_j=\Theta(\la (n-\ell))$; and for each $i=2,3,4$, $\sum_{\ell<j\le n}[d_j]_i=O(\la^i (n-\ell))$. Many degree sequences satisfy such conditions including interesting examples in which
 the last $n-\ell$ vertices have the same degree; or their degrees are highly concentrated; or the degree distribution is Poisson-like or truncated-Poisson-like. Then, with basically the same proof as Theorem~\ref{thm:bidegrees}, we will have
\[
g({\bf d})=\frac{|\Phi|}{\prod_{i=1}^n d_i!}\exp\left(-\frac{M_1}{2}+\frac{M_2}{2M_1}+\frac{3}{4}+\sum_{i<j} \log(1+d_id_j/M_1) +O\Big(\sqrt{\xi}\Big)\right),
\]
if Theorem~\ref{thm:bidegrees} (a) holds with $\delta$ replaced by $\la$ and $\Delta\ell=O(\la n)$; or if Theorem~\ref{thm:bidegrees} (b) holds with $\delta$ replaced by $\la$.
 \medskip

 In a typical power-law  density-bounded sequence (i.e.\ with $2<\gamma<3$), the maximum degree is $o(M_1^{1/2})$, as discussed above, whilst for  a power-law distribution-bounded sequence,  it can reach $M_1^{1/\sqrt 3}\approx M_1^{0.577}$. To illustrate the power of our main result, we consider a  generalised type of power-law degree sequence  with an even stronger  tail in its distribution, in which  the maximum degree can approach $M_1^{3/5}$.

 We say ${\bf d}=(d_1,\ldots, d_{n})$ follows a {\em long-tailed power law} with parameters $(\alpha,\beta,\gamma)$, if there is a constant $C>0$ such that for every $n$, 
\begin{enumerate}
\item[(a)] each coordinate is non-zero and either  at most $C$  or at least $n^{\alpha}$;
\item[(b)] for every integer $i\ge 1$, the number of coordinates whose value is at least $i n^{\alpha}$ but less than $(i+1) n^{\alpha}$ is  at most  $Cn^{\beta} i^{-\gamma}$.
\end{enumerate}
 Note that when $\alpha=0$ and $\beta=1$, this definition agrees with that of power-law  density-bounded  degree sequences.

\begin{thm}\lab{thm:heavyPowerLaw}
Let ${\bf d}$ be a long-tailed power-law degree sequence with parameters $(n,\alpha,\beta,\gamma)$ such that $1<\gamma<3$,  $\gamma\ne 2$,   $\alpha>1/2$ and
$$
 0< \beta<\left\{
\begin{array}{rl}
\displaystyle\frac{3-5\alpha}{1+6/\gamma}& \mbox{if\/ $2\le \gamma<3$}\\
\raisebox{6mm}{\phantom{z}}\displaystyle\frac{3-5\alpha}{8/\gamma}& \mbox{if\/  $1<\gamma<2$}.
\end{array}
\right.
$$
Then
\[
g({\bf d})=\frac{|\Phi|}{\prod_{i=1}^n d_i!}\exp\left(-\frac{M_1}{2}+\frac{M_2}{2M_1}-\frac{M_3}{3M_1^2}+\frac{3}{4}+\sum_{i<j}\Big(\log(1+d_id_j/M_1)\Big)+O\Big(\sqrt{\xi}\Big)\right),
\]
where
\[
  \xi=\left\{
  \begin{array}{ll}
  n^{5\alpha+\beta+6\beta/\gamma-3}   & \mbox{if $2< \gamma<3$}\\
   n^{5\alpha+4\beta -3}\log n   & \mbox{if $\gamma=2$}\\
   n^{5\alpha+8\beta/\gamma-3}  & \mbox{if $1<\gamma<2$},
  \end{array}
  \right.
  \]
which is $o(1)$ by the assumption on $\beta$.
 \end{thm}
\no {\bf Remark}. For convenience, we omitted the  case $\gamma\ge 3$ which can be easily worked through if required. We also omitted the case   $\alpha<1/2$, even though the forthcoming main result will still apply under appropriate conditions, because those conditions are much more complicated.

Next we  present our general result, from which the foregoing special cases are derived.   Define
\bea
U_1&=&\sum_{v\le n}(d_v-2)\min\{[d_v]_2/M_1,1\};\non\\
U_2&=&\sum_{1\le u<v\le n }\min\{[d_u]_2[d_v]_2/M_1^2,d_ud_v/M_1\};\non\\
U_3 &=& \sum_{u\ne v\le n}\sum_{w \le n}
\min\{[d_u ]_2[d_v]_2/M_1^2,d_ud_v/M_1\}\min\{[d_u-2]_2[d_v]_2/M_1^2,1\}(d_w-2);\lab{tildeUs} \\
U_4 &=& \sum_{u\ne v\le n}  \min\{[d_u]_3[d_v]_2/M_1^2,[d_u]_2d_v/M_1\};\non\\
U_5 &=& \sum_{u\ne v\le n}\sum_{w \le n}
\min\{ d_u [ d_v ]_2/M_1^2, d_v  /M_1\} \min\{ [d_u-2]_2[d_w]_2/M_1^2,(d_u-2)d_w/M_1\}.\non
\eea
\begin{thm}\lab{thm:general}   Let $U_k$  be defined as above for $1\le k\le 5$. Define
\[
 \xi= U_5+
\frac{U_1+ U_2^2+ U_3}{M_1 }
+\frac{     U_4  M_2}{M_1^2}
+\frac{   U_2  M_2^2}{M_1^3 }
+\frac{ M_2}{M_1^2}+\frac{  M_3  M_2 }{M_1^3}+\frac{ M_2^3}{M_1^4}.
\]
 Suppose that $\xi=o(1)$. Then
\bean
g({\bf d})&=&(1+O(\sqrt{\xi}+M_1^{-1}))\frac{|\Phi|}{\prod_{i=1}^n d_i!}\exp\left(-\frac{M_1}{2}+\frac{M_2}{2M_1}-\frac{M_3}{3M_1^2}+\frac{3}{4}+\sum_{i<j} \log(1+d_id_j/M_1) \right).
\eean
\end{thm}

 Since $\xi$ has a rather complicated formula, we give a simple upper bound on $\xi$ next.

\begin{lemma}\lab{lem:xisimple1}
Let $\xi$ be defined as in Theorem~\ref{thm:general}. Then
\begin{enumerate}
\item[(a)] $\xi = O\Big( (M_2+ M_3)/M_1^2+(M_2^2 M_3 +M_2^3)/M_1^4+(M_2^4+ M_2 M_3 M_4)/M_1^5\Big);$
\item[(b)] if $\Delta=O(\sqrt{M_1})$, then the term $M_2 M_3 M_4/M_1^5$ in (a) can be dropped, and the resulting bound on $\xi$ is tight to within a constant factor.
\end{enumerate}
\end{lemma}

We can often get better bounds on $\xi$ when additional constraints are placed on  the degree sequence (particularly when $\Delta=\Omega(\sqrt{M_1})$). These bounds will be presented in Section~\ref{sec:boundxi}.
In the next section, we prove Theorem~\ref{thm:general}. We will derive Theorems~\ref{thm:M2}--\ref{thm:heavyPowerLaw} as special cases of Theorem~\ref{thm:general}  in Section~\ref{sec:applications}.

 We close this section with a short verification that the main term of Theorem~\ref{thm:M2} agrees with~\cite[Theorem 1.4]{J} for  $M_2=O(M_1)$. The latter result gives  an asymptotic formula for $g({\bf d})$ in which the logarithm of $\pr(\G^*(n,{\bf d})\ \mbox{is simple})$  is expressed  as
\bel{janson}
-\frac12 \sum_i\mu_{ii} - \sum_{i<j}\big(\mu_{ij}-\log(1+\mu_{ij})\big)
\ee
where
$$
\mu_{ii}=[d_i]_2/M_1,\quad  \mu_{ij}=\sqrt{[d_i]_2[d_j]_2}/M_1.
$$
To see that this is within $o(1)$ of the exponent in~\eqn{g(d)}  is not entirely straightforward, so we give some details of an argument.
Let 
$ 
F(\la)=\la -\log(1+\la) 
$ 
and write $\la$ for $\la_{ij}$, and similarly $\mu$. Consider two cases depending on some fixed $0<\eps<1/12$.
Firstly, if $d_id_j\le M_1^{1-\eps}$, then $\la$ and $\mu$ are  both $O(M_1^{-\eps})$, and subtracting the resulting expansions gives 
\bean
F(\mu)-F(\la)= -\frac{d_id_j(d_i+d_j-1)}{2M_1^2} + O\big(d_i^2d_j^2(d_i+d_j)/M_1^3+M_1^{-3}\big).
\eean
Secondly, if
 $d_id_j>M_1^{1-\eps}$, then since $d_i$ and $d_j$ are both at most $\Delta =O(\sqrt{M_2}) =O( M_1^{1/2})$, they are both  $\Omega(M_1^{1/2- \eps })$, and thus 
  $\mu = \la \big(1+O(d_i^{-1}+d_j^{-1})\big) =\la+O(M_1^{ -1/2})$.  Now $F'(\la)=O(1)$, so  
$
 F(\mu)-F(\la)=O(M_1^{ -1/2}).
 $ 
 Since $M_2=O(M_1)$ and $d_i(d_i-1)=   \Omega(M_1^{1- 2\eps })$, there are only $O(M_1^{4\eps})$  values of $(i,j)$ in this case, giving a total error of $O(M_1^{4\eps})O(M_1^{-1/2})=o(1)$. So we can use $F(\la)$ in place of $F(\mu)$ for these terms.  Moreover, the sum of $d_id_j(d_i+d_j-1)/2M_1^2$, over these terms, is $o(1)$, so we can combine the two cases to obtain 
$$
 \sum_{i<j}\big(\mu_{ij}-\log(1+\mu_{ij})\big)=o(1)+ \sum_{i<j}\big(\la_{ij}-\log(1+\la_{ij})\big)- \sum_{i<j} \frac{d_id_j(d_i+d_j-1)}{2M_1^2}.
 $$
The last term here is $o(1)-M_1^{-2} \sum_{i,j}(  d_i^2d_j/2-d_id_j/4)$ (noting that terms with $i=j$ are $O(M_3/M_2^2)=o(1)$). We have
$$
\sum_{i,j}(  d_i^2d_j/2-d_id_j/4)=\sum_{i,j}\big(d_i(d_i-1) d_j/2+d_id_j/4\big)
= \frac12 M_2/M_1+\frac14.
$$
The first term in~\eqn{janson} is 
$-\frac12  \sum_i\mu_{ii}=-\frac12 M_2/M_1$. 
Hence~\eqn{janson} is  
$
o(1)+  \frac14 - \sum_{i<j}\big(\la_{ij}-\log(1+\la_{ij})\big).
$
Up to $o(1)$ terms, this agrees with~\eqn{g(d)}, since   $M_3/M_1^2=o(1)$ and, as noted above, $-\sum_{i<j}\la_{ij}= \frac12(M_2/M_1-M_1+1)$. 
 
 \section{Proof of Theorem~\ref{thm:general}}
\lab{sec:mainproof}

 As we assume all $d_i$ are at least 1, we may assume without loss of generality  that  $d_1\ge \cdots \ge d_n \ge 1$. Recall
that $\Phi$ denotes the set of pairings with   degree sequence ${\bf d}$.  Let $\P\in\Phi$. We often refer to the multigraph corresponding to $\P$ as if it were the same as $\P$, and hence we sometimes call  the bins in $\P$ vertices, and treat the pairs in $\P$ as edges.  For two (possibly equal) vertices $u$ and $v$, we say that $uv$ is a {\em multiple edge}, of multiplicity $i$, if there are $i\ge 2$ pairs with end-vertices $u$ and $v$.
A {\em single} edge is an edge that is not (part of) a multiple edge, and a {\em loop} has both ends at the same vertex.

Let
${\nat}_{\ge k}$ denote  the set of integers at least $k$. In this paper, we use matrices whose entries are not just numbers, but can be $\bit$ as well.  Define  $\M $ to be the set of $n\times n$ symmetric matrices  ${\bf M}=(m_{i,j})$ for which $m_{i,j}\in\{\bit\}\cup {\nat}_{\ge 2}$ if $i<j$ and    $m_{i,j}\in {\nat}_{\ge 0}$  if $i=j$.

Given $\P\in\Phi$, we define the {\em signature} of $P$ to be the matrix $\MM(P)\in\M$ defined as follows. For $i\ne j$, if the multiplicity of the edge $ij$ in $P$ is at least 2, then $m_{i,j}$ equals that multiplicity, whilst if the multiplicity is 1 or 0, $m_{i,j}= \bit$. For each $i$, $m_{i,i}$ is the number of loops at $i$. Next, for any  ${\bf M}\in \M$, let $\C({\bf M})$ be the set of $\P\in\Phi$ whose signature is $\MM$. Then for any $\P\in\C(\MM)$, the locations and multiplicities
 of all loops and multiple edges in $G(\P)$ are determined. Note that the single non-loop edges in this graph are unconstrained apart from the number  of such edges  incident with each vertex.

 Define $\MMS$ to be the matrix in $\M$  with $\bit$ in all off-diagonal
 positions and $0$ on the diagonal. Thus, a pairing is simple if and only if its signature is $\MMS$.

 Note that $\pr(\G^*(n,{\bf d})\ \mbox{is simple})=|\C(\MMS)|/|\Phi|$. We will estimate this ratio using switchings. In this argument, we often need a bound of the number of 2-paths starting from a given vertex $v$ in a pairing, in order to bound the number of possible switchings from below. The trivial upper bound for  the number of  2-paths is $\Delta(\Delta-1)$, which would place a natural restriction on $\Delta$ as in many previous works (see for instance~\cite{MW2}). Since at most $d_i-1$ of the 2-paths use vertex $i\ne v$, we can use another simple (and clearly not sharp) upper bound: $\tau:=\sum_{i=1}^{\Delta} d_i$.
  Before proceeding, it is useful and informative to obtain bounds on  $M_2$,  $\Delta$ and  $\tau$ in terms of $M_1$ in the setting of  Theorem~\ref{thm:general}.
\begin{lemma}\lab{lem:tau}
 Define $
\xi$ as in Theorem~\ref{thm:general}. If $\xi=o(1)$ then $M_2=o(M_1^{4/3})$, $\Delta=o(M_1^{3/5})$ and  $\tau=o(M_1)$.
\end{lemma}

\proof  The first and second bounds come immediately from considering the terms   $M_2^3/M_1^4$ and $M_3M_2/M_1^3$ in $\xi$.

Let $c>0$ be an arbitrary constant. Suppose to the contrary that $\tau>c M_1$. Then,  using Cauchy's inequality and also $\Delta=o(M_1^{3/5})$, we have
$$
M_2\ge \sum_{i=1}^{\Delta}d_i^2 - M_1 \ge \frac{\tau^2}{\Delta}-M_1=\Omega(M_1^{4/3}).
$$
However, this contradicts the fact that $M_2=o(M_1^{4/3})$. Hence, $\tau=o(M_1)$. \qed

\subsection{Auxiliary functions} \lab{sec:auxiliary}
 Recall from~Section~\ref{sec:mainproof} that given a symmetric matrix  ${\bf M}=(m_{i,j})\in \M$, $\C({\bf M})$ is the set of $\P\in\Phi$ whose signature is $\MM$.  Note for later use that
\bel{Mtophi}
\bigcup_{{\bf M}\in\M }  \C({\bf M})  =  |\Phi|.
\ee
We will  use the following auxiliary functions of ${\bf M}\in \M$,  some of which also  depend on a pair $(i,j)$ which in our applications will be two distinct vertices.
\begin{itemize}
\item  $Z ({\bf M})=\sum_{1\le u<v\le n }I_{m_{u,v}\ge 2}\, m_{u,v}$, where $I_{B}$ is the characteristic function or indicator of an event $B$. This is the total number of pairs in  non-loop  multiple edges  for pairings in $\C ({\bf M})$.
 \item $Z_{i,j} ({\bf M}) = Z ({\bf M})-m_{i,j}$.
\item  $Z_2({\bf M}) =\sum_{1\le u<v\le n  }I_{m_{u,v}\ge 2}\,m_{u,v}^2$.
\item $Z_0({\bf M})=\sum_{1\le u\le n } m_{u,u}$.
 \item $W_{i,j}({\bf M})=\sum_{w }(d_w-2)$   where the summation is over all  $w\ne j$ such that  $m_{i,w}\ge 2$, which effectively means that $iw$ is  designated as a multiple edge   by ${\bf M}$.
 \item $Q_{i,j}({\bf M})=\sum_{\{u,v\}}(d_u-2)(d_v-2)$ where the summation is over all pairs $\{u,v\}\ne \{i,j\}$ such that $u<v$ and $m_{u,v}\ge 2$. (We exclude $\{i,j\}$  because our argument later requires the set of matrix entries that influence $Q_{i,j}$ to be independent of $m_{i,j}$.) 
\item $R_{i,j}({\bf M})=\sum_{w} m_{i,w}$  where the summation is over all $w$ with $ m_{i,w} \ge 2$ and $w\notin\{i, j\}$.
\end{itemize}

When convenient, we abbreviate $Z ({\bf M})$ to $Z $, and similarly for the other variables defined above.  Note that all functions $Z_{i,j}$, $W_{i,j}$, $Q_{i,j}$ and $R_{i,j}$ are independent of $m_{i,j}$.  We will use this property later in our argument. Indeed, this is the motivation for the definition of both  $Z_{i,j}$ and $Z$.

\subsection{Switchings for multiple edges}
\lab{sec:multiple}
  In this subsection, we deal with multiple non-loop edges. Consider the following assumptions on ${\bf M}$, where $\xi_1$ is a certain function that is $o(1)$, to be specified later.
\begin{enumerate}
\item[(A1)] For every $i< j$, $m_{i,j}^2\le  \xi_1M_1$.
 \item[(A2)] $Z({\bf M})\le\xi_1M_1$.
 \end{enumerate}
One more assumption (A3) is to be stated in the proof. We will later show that,  for a random $\P\in \Phi$, the probability that   ${\bf M}(\P)$ fails any of these assumptions tends to 0 quickly. We now fix a matrix ${\bf M}=(m_{i,j})\in \M$ satisfying the three assumptions.

Next fix a pair $(i,j)$  with $i<j$ for which $m_{i,j}\ge 2$  (and thus $m_{j,i}\ge 2$ since $\MM$ is symmetric).
 For $m\ge 0$,  let
${\bf M}(m)$   be  the  matrix which is formed by changing the  $(i,j)$ and $(j,i)$ entries  of  ${\bf M} $ from $m_{i,j}$ to $m$.
 We extend the definition of $\C(\MM(m))$ in the obvious way to the cases $m=0$ and 1 (where $\MM(m)\notin \M$).   Similarly define ${\bf M}(\bit)$.  Note that $\C(\MM(\bit))=\C(\MM(0))\cup\C(\MM(1))$ and $\MM(\bit)\in\M$.   Note  also  that  ${\bf M}(m)$  depends on the values of $i$ and $j$, but they are fixed so we suppress them from the  notation.
 We first use the switching method to estimate the ratio
$\big|\C\big({\bf M}(m)\big)\big|/\big|\C\big(  {\bf M}(0)\big) \big|$.

 Let $m\ge 1$ with $m\le m_{i,j}$, so that by (A1) we have $m^2=o(M_1)$.
 For any $\P\in \C\big({\bf M}(m)\big)$, define a ``switching'' operation as follows.
Label the endpoints of the $m$ pairs between $i$ and $j$ as $2g-1$ and $2g$, $1\le g\le m$, where points $1,3,\ldots,2m-1$ are contained in vertex $i$.
Pick another $m$ distinct pairs $x_1,\ldots,x_m$. Label the endpoints of
$x_g$ as $2m+2g-1$ and $2m+2g$. The switching operation replaces pairs $\{2g-1,2g\}$ and $\{2m+2g-1,2m+2g\}$ by $\{2g-1,2m+2g-1\}$ and $\{2g,2m+2g\}$.
See Figure~\ref{f:heavy} for an example when $m =2$.

\begin{figure}[htb]
\vbox{\vskip -.8cm
 \hbox{\centerline{\includegraphics[width=12cm]{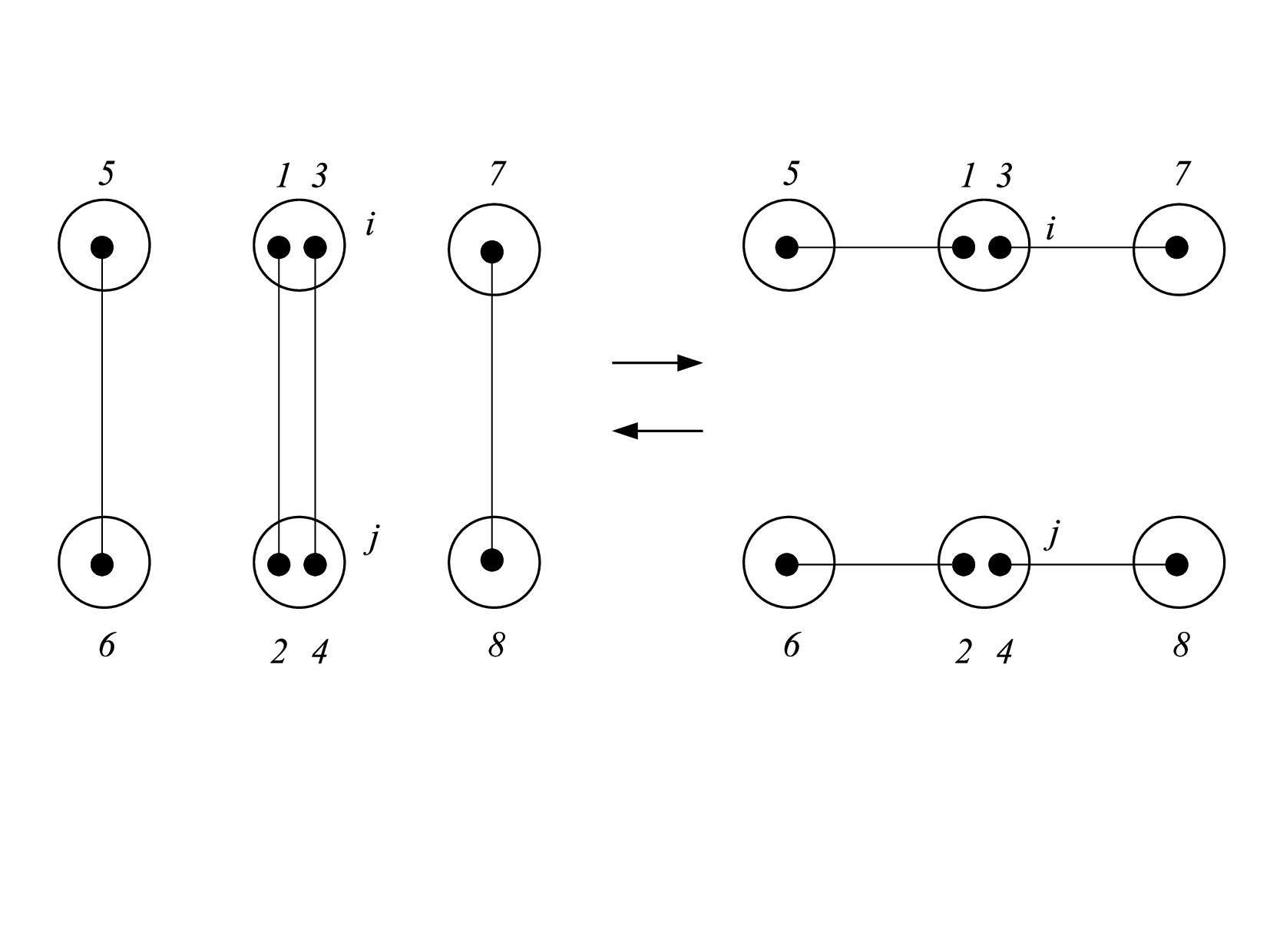}}}

\vskip -2.5cm \smallskip} \caption{\it  switching for multiple edges}

\lab{f:heavy}

\end{figure}

   It is easy to see that this switching converts the pairing $\P\in\C\big({\bf M}(m)\big)$
to a pairing in  $\C\big({\bf M}(0)\big)$, provided that none of the following conditions hold  (though these are not all entirely necessary):
\begin{description}
\item {(i)}  a pair $x_g$ is part of a non-loop multiple edge;
\item {(ii)} a pair $x_g$ uses a vertex $w$ already adjacent in $\P$  to $i$ or $j$ by a multiple edge, and $x_g$ does not satisfy (i) (this would increase the multiplicity of a multiple edge);
\item {(iii)} a pair $x_g$ uses a vertex $w$ already adjacent in $\P$  to $i$ or $j$  by a single edge (this would create a new multiple edge);
\item{(iv)} some two pairs $x_g$ and $x_{g'}$ have a common end-vertex (if this were permitted, two of the new pairs can possibly create a multiple edge).
 \item{(v)} a pair $x_g$ is incident with $i$ or $j$ or is a loop.
\end{description}

 We call a switching satisfying these conditions {\em good}.

We will bound the probability that a randomly chosen switching is not good, when applied to a {\em random}  $\P\in\C\big({\bf M}(m)\big)$. In all cases but (iii)  our bound actually applies to  an arbitrary $\P$ rather than a random pairing.

   We can choose the pairs one by one. Each of them is potentially any one of the $M_1/2$ possible pairs, but it must avoid the $m$ pairs joining $i$ and $j$, as well as up to $m-1$ other pairs already chosen, so  there are $M_1/2-O(m)$ options.
Since $m=o(M_1)$, the probability that a randomly chosen $x_g$ is any given pair, conditional on the previous pairs $x_1,\ldots, x_{g-1}$, is $O(1/M_1)$.
  We will use this observation several times. In particular,   for each $1\le g\le m$, conditional on the choices of $x_1,\ldots, x_{g-1}$, the probability that $x_g$ is one of $x_i$, $i\le g-1$, or $x_g$
   is a pair between $i$ and $j$, is $O(m/M_1)$. Taking the union bound over all $1\le g\le m$, the probability that $x_g$'s are not distinct, or use a pair between $i$ and $j$, is $O(m^2/M_1)$. By the above observation, this probability is $o(1)$.

   Since at most $Z_{i,j}$
   pairs are in multiple edges of ${\bf M}(0)$, the probability of (i) occurring is  $O(mZ_{i,j} /M_1)$.

The number of pairs that cause the condition (ii) for any $x_g$ is the sum of $d_w-2$ over all $w$ such that $iw$ or $jw$ is a multiple edge (excluding $w=i$ or $j$), which is $W_{i,j}+W_{j,i}$. Arguing as for (i), the probability of this occurring is  at most $m(W_{i,j}+W_{j,i}) /M_1$.

For condition (iii), we need to argue about the expected number of switchings in which the condition occurs, when applied to a random pairing $\P\in\C\big({\bf M}(m)\big)$. To this end, we first estimate the probability $p(i,w)$ that a given point $a$ in vertex $i$ is paired in $\P$ to a given point $b$ in vertex $w$.
This uses the following subsidiary switching argument.

For those pairings $\P$ containing the pair $ab$, consider switching $ab$ with another randomly chosen pair $a'b'$ in $\P$, i.e.\ delete the pairs $ab$ and $a'b'$ and insert the pairs $aa'$ and $bb'$ to create a new pairing $\P'$.  Then $\P'$ is also
in $\C\big({\bf M}(m)\big)$ provided that neither $a'$ nor $b'$ is in a vertex adjacent to $i$ or $j$ (there are at most $2\tau$ such points since each corresponds to a unique 2-path starting from $i$ or $j$) and $a'b'$ is not in a  multiple edge (there are at most $2Z$ such points).  By assumption (A2) and  Lemma~\ref{lem:tau},
 the number of ways to choose $a'b'$ such that $\P'$ is also
in $\C\big({\bf M}(m)\big)$ is  $\Omega(M_1)$.
Furthermore, each such $\P'$ can be created in at most one way.  Hence $p(i,w)=O(1/M_1)$.

Applying the union bound to all appropriate  $(a,b)$ (where we can assume $a$ is not in one of the pairs joining $i$ and $j$), we find that the probability that $iw$ is  an edge in $\P$ is $O\big((d_i-2)d_w/M_1)$. Conditional upon this event, the probability that the  random pair $x_g$ chooses a point in $w$ is   $O( d_w/M_1)$.
   The same considerations as for $i$ apply to $j$, and we conclude that  the probability of (iii) occurring for a random $\P$ is $O(m\sum_{w\le n}(d_i-2+d_j-2) d_w^2/M_1^2)=   O\big(m(d_i+d_j-4) (M_2+M_1)/M_1^2\big)$.

Given  an arbitrary  $\P$,  a randomly chosen switching satisfies (iv) with probability $O(m^2  M_2/M_1 ^2) $ since there are $O(m^2)$ choices of $(g,g')$,   and the number of ways that any two of  them
can both choose a point in the same vertex is   $O( M_2)$.

 Finally, for an arbitrary $\P$, a randomly chosen switching satisfies (v) with probability $O\big(m(d_i+d_j-4+Z_0)/M_1\big)$, since for each $g$, the probability that $x_g$ contains a point in $i$ or $j$ is $O((d_i+d_j-4)/M_1)$ and the probability it is a loop is $O( Z_0 /M_1)$. Here the term $-4$ occurs because there are at least four points in the multiple edge joining $i$ and $j$ that are excluded.

Let $N_{i,j}$ denote the expected number of choices of  a good  switching, for a  uniformly randomly chosen $\P\in\C\big({\bf M}(m)\big)$, where we distinguish between the $m!$ different ways to assign the labels $1,3, \ldots,2m-1$  to the points in vertex $i$. These induce labels of the points paired with them in $j$. There are two ways to label the ends of each of the chosen pairs, and, using the observation before considering (i), there are  $\big(1-O(m^2/M_1)\big)M_1^m/2^m$ ways to choose the pairs. Hence
\bea
N_{i,j}&=& m! M_1  ^m  \bigg(1+O(m)\bigg( \frac{Z_{i,j}}{M_1}+  \frac{W_{i,j}+W_{j,i}}{M_1}  +   \frac{(d_i+d_j-4) M_2}{M_1^2}+\frac{d_i+d_j-4+Z_0}{ M_1}\bigg)  \non \\
&&+O(m^2)\bigg(   \frac{ M_2 }{M_1 ^2}+   \frac{1 }{M_1} \bigg)\bigg).\lab{Nij}
\eea

On the other hand, speaking informally,  the  inverse of a good switching will convert a pairing   $\P\in \C\big(  {\bf M}(0) \big)$ to one in  $\C({\bf M}(m))$.  We formally define an inverse switching to be the following operation. Pick $m$ distinct
points in vertex $i$  and label them as $ 2g-1$, $g=1,\ldots,m$; label the point paired with $ 2g-1$ as   $2m+2g-1$;
 do a similar thing for $j$, producing pairs $\{2g,2m+2g\}$ as in the right hand side of Figure~\ref{f:heavy}; and finally
replace the pairs $\{2g-1,2m+2g-1\}$ and $\{2g,2m+2g\}$ by new pairs $\{2g-1,2g\}$ and $\{2m+2g-1,2m+2g\}$ for all appropriate $g$.
An inverse switching is called good if it creates a pairing in $\C\big(  {\bf M}(m) \big)$, i.e.\ if it has the reverse effect of a good switching applied to a pairing in $\C\big(  {\bf M}(m) \big)$. Since no pair in $\P$ joins vertices $i$ and $j$, the inverse switching is good if none of the following conditions holds:
\begin{description}
\item {(vi)}  a pair picked incident with $i$ or $j$ is part of an existing multiple edge or forms a loop;
\item {(vii)} a new pair is added in parallel to an existing multiple edge.
\item {(viii)} a new pair is added in parallel to an existing single edge.
\item {(ix)} a pair picked incident with $i$ has a common end vertex with a pair picked incident with $j$. (Then a new loop would be created.)
\end{description}

We next bound the probability that a randomly chosen inverse switching is not good, when applied to a  random   $\P\in\C\big({\bf M}(0)\big)$. Note that there are potentially $[d_i]_m[d_j]_m$   inverse switchings, but some of these may not be good.

 First consider (vi).  The number of   pairs incident with vertex $i$ that are already part of a multiple edge is  $R_{i,j}$. For those already part of a loop, it is $m_{i,i}$. Thus, the proportion of the initial count of switchings falling into this case is
   $O(m)\big((R_{i,j}+m_{i,i})/d_i + (R_{j,i}+m_{j,j})/d_j\big)$.

For (vii), again we consider a random pairing.  Suppose the points $2g-1$ in vertex $i$ and  $2g$  in vertex $j$ are specified,
and consider the random pairing $\P\in\C\big({\bf M}(0)\big)$. Given a multiple edge $uv$, the probability that $2g-1$ and $2g$ are paired with points in $u$ and $v$ respectively is, arguing as for (iii) and switching out the two pairs simultaneously, $O((d_u-2)(d_v-2)/M_1^2)$. (Here we use $d_u-2$ rather than $d_u$ since the points paired with $2g-1$ and $2g$ cannot be part of the multiple edge, which has multiplicity at least 2.) Hence, the probability of (vii) is $O(mQ_{i,j}/M_1^2)$.

 For (viii),   we do not need to consider cases which also fall into (vi). Again we   need to consider a random pairing.
 Suppose the points $2g-1$ in vertex $i$ and  $2g$  in vertex $j$ are specified, and consider the random pairing $\P\in\C\big({\bf M}(0)\big)$.  Conditioning on the two pairs containing these points, the rest of $\P$ is a uniformly random pairing conditional upon all the multiple edges existing as specified by ${\bf M}(0)$. Let $u$ and $v$ be two vertices.   The number of ways to choose an ordered pair of points $(a,a')$ and $(b,b')$ in each of $u$ and $v$  respectively is $[d_u]_2[d_v]_2$. Arguing as for (iii), switching three pairs out at once, the probability that $2g-1$ and $2g$ are paired with $a$ and $b$, and $a'$ and $b'$ form a pair, is $O(1/M_1^3)$.   Summing over all $u$ and $v$, we obtain the bound $ M_2^2/M_1^3$ on the probability that any two chosen points $2g-1$ and $2g$ lead to (viii).  By the union bound, the probability that at least one of the $m$ new edges causes condition (viii) is   $O(m{ M_2}^2/M_1^3)$.

For (ix), we argue as for (vii), but noting that the probability that  $2g-1$ and $2g$ are paired with points in the same vertex $w$ is $O\big(d_w(d_w-1)/M_1^2\big)$. Thus, the bound for (ix) is $O(m{M_2} /M_1^2)$.

 Letting $N'_{i,j}$ denote the expected number of good inverse switchings for a random $\P\in\C\big({\bf M}(0)\big)$,  we get
\bel{Nijprime}
N'_{i,j} = [d_i]_m[d_j]_m \bigg(1+O(m)\bigg( \frac{R_{ i,j}+m_{i,i} }{d_i} +\frac{R_{j,i}+m_{j,j} }{d_j}      +  \frac{Q_{i,j}}{M_1^2} +   \frac{{M_2}^2}{M_1^3} + \frac{{M_2}}{M_1^2} \bigg)\bigg).
\ee
\subsection{Eliminating multiple edges}
Before proceeding, we let  $\eta_{i,j }({\bf M},m)$ denote the sum of the error terms in~\eqn{Nij} and~\eqn{Nijprime}, i.e.\ $\eta_{i,j }({\bf M},m)= mZ_{i,j}/M_1+\cdots+ m{M_2} /M_1^2 $.
The last assumption that we will make on the matrix ${\bf M}$ to which we apply the switching analysis is the following.
 \begin{enumerate}
\item[(A3)] \ $\eta_{i,j}(\MM,m_{i,j})\le \xi_1$ for every $ i<j $  such that $m_{i,j}\ge 2$.
\end{enumerate}

 Moreover, we apply the switchings only in the case that $m=m_{i,j}$ or $m=1$.   As with earlier notation,  we use $\eta_{i,j }(m)$ to denote $\eta_{i,j }({\bf M},m)$. For each $m\ge 1$,   if we let
$$
\rho_m = \frac{\big|\C\big({\bf M}(m)\big)\big|}{\big|\C\big({\bf M}(0)\big)\big|}, 
$$
we have from~\eqn{Nij} and~\eqn{Nijprime}
\bel{rho}
\rho_m = \frac{[d_i]_m[d_j]_m}{m!M_1^m}\big(1+O(\eta_{i,j }(m))\big).
\ee

Next, write $C_m$ for $\big|\C\big({\bf M}(m)\big)\big|$ (where $m\ge 0$).
Then
\bea
\frac{\big|\C\big({\bf M}\big)\big|}{\big|\C\big({\bf M}(\bit)\big)\big|}&=&\frac{C_{m_{i,j}}}{C_0+C_1} =  \frac{ \rho_{m_{i,j}}}{ 1+  \rho_1}\non\\
 &=&  \frac{[d_i]_{m_{i,j}}[d_j]_{m_{i,j}} }{(1+ d_id_j /M_1){m_{i,j}}!M_1^{m_{i,j}}} \big(1+ O(\eta_{i,j }({m_{i,j}}) )\big) \lab{ratclub}
 \eea
since
  $$
 1+\rho_1= 1+ d_id_j / M_1 +O(\eta_{i,j }(1)d_id_j / M_1) =(1+ d_id_j / M_1)(1+ O(\eta_{i,j }(1)d_id_j /(M_1+d_id_j) ),
  $$
and it is easy to check that $\eta_{i,j }(1)\le \eta_{i,j }( m_{i,j})$  since  ${m_{i,j}}\ge 2$ and the functions $Z_{i,j}$, $W_{i,j}$ etc.\   are functions of $\MM$,  independent of $ m_{i,j} $.

 Recall that $(i,j)$ is fixed.  The equation above estimates the effect on $\big|\C\big({\bf M} \big)\big|$ of changing an $(i,j)$-entry of $\bf M$ (and simultaneously $(j,i)$, to keep $\MM$ symmetric) from a number at least 2, to  $\bit$. Next, we can select another non-$\bit$ entry of ${\bf M}$ and change it (and the symmetric entry)  to $\bit$ using the same procedure.  Let ${\bf M} _\bit$ be the matrix with all off-diagonal entries equal to $\bit$  and  each  $(i,i)$ entry equal to $m_{i,i}$, and let $\Heavy({\bf M})= \{(i,j):   i<j, m_{i,j}\ge 2\}$. Then, applying~\eqn{ratclub} for each $(i,j)\in \Heavy({\bf M})$, we obtain the    formula
\bel{formula}
 \frac{\big|\C\big({\bf M} \big)\big|}{\big |\C\big({\bf M} _\bit\big)\big|}
 =   \exp\big( O\big(\eta({\bf M} )\big)\big)\prod_{(i,j)\in \Heavy({\bf M})} \frac{[d_i]_{m_{i,j}}[d_j]_{m_{i,j}} /(m_{i,j}!M_1^{m_{i,j}})}{1+ d_id_j /M_1},  \ee
as long as
 $$
 \sum_{(i,j)\in \Heavy({\bf M})}\eta_{i,j }(m_{i,j}) = O\big(\eta({\bf M})\big).
 $$
 Note that since  $\sum_{(i,j)\in \Heavy({\bf M})}m_{i,j } = Z({\bf M})$, $Z\le Z_2$,  and since $d_i\ge d_j$ as $i<j$, terms in  $\sum_{(i,j)\in \Heavy(\MM)} \eta_{i,j}(m_{i,j})$  like $d_j/M_1$ drop out,  and  we can take
\bea
&&\hspace{-0.8cm}\eta\big({\bf M} \big) =
  \frac{  Z_2+ZZ_0}{M_1 }  +\frac{M_2Z_2}{M_1 ^2}+\frac{ M_2 ^2Z}{M_1^3}   \non\\
  &&\hspace{-0.4cm}+
\sum_{(i,j)\in \Heavy({\bf M}) } m_{i,j}\Bigg(
\frac{Z_{i,j}+  W_{i,j} + W_{j,i} +d_i-2 }{M_1 }+ \frac{(d_i -2) M_2 }{M_1^2 }    +\frac{R_{i,j}+m_{i,i}}{d_i}+\frac{R_{ j,i}+m_{j,j}}{d_j} + \frac{Q_{i,j}}{M_1^2} \Bigg).\lab{eta}
\eea

\subsection{Switchings for loops}

 Recall \label{loopswitching} that given $\MM\in \M$, $\MM_{\bit}$ is the matrix by changing all  off-diagonal   entries in $\MM$ to $\bit$.
Next, we estimate the ratio $|\C(\MM_\bit)|/|\C(\MMS)|$, where $\MMS$ as defined earlier is the  symmetric matrix with $\bit$ in all off-diagonal
 positions and $0$ on the diagonal.

 Fix $1\le i\le n$ with $m_{i,i}\ge 1$.
Let $\MM(0)$ be the matrix obtained from $\MM_\bit$ by changing the $(i,i)$ entry from $m_{i,i}$ to $0$. We use $m$ to denote $m_{i,i}$ for convenience.
To eliminate the loops at vertex $i$ in $\C(\MM_\bit)$, we define the following switching. Take $\P\in \C(\MM_\bit)$ and label the $m$ loops at $i$ by $1,\ldots,m$. For the $g$th loop at $i$, label the endpoints   $2g-1$ and $2g$ (for $1\le g\le m$). Pick $m$ distinct pairs in $\P$,  labelling the endpoints of the $g$th pair  $2m+2g-1$ and $2m+2g$,  and then pick another $m$ distinct pairs and label the endpoints of the $g$th of this lot of pairs $4m+2g-1$ and $4m+2g$.  The switching replaces pairs
$\{2g-1,2g\}$, $\{2m+2g-1,2m+2g\}$ and $\{4m+2g-1,4m+2g\}$ by $\{2g-1,2m+2g-1\}$, $\{2g,4m+2g\}$ and $\{2m+2g+4m+2g-1\}$. See Figure~\ref{f:loop} for an illustration. Let $v_j$ denote the vertex containing point $j$ for all $2m+1\le j\le 6m$.
 We say a switching is good if none of the following conditions holds:
 \begin{enumerate}
\item[(a)] $v_j= i$ for some $2m+1\le j\le 6m$, or the vertices in the set $\{v_j,\ 2m+1\le j\le 6m\}$ are not all distinct;
\item[(b)] $i$ is already adjacent to some $v_{2m+2g-1}$ or to some $v_{4m+2g}$;
\item[(c)] $v_{2m+2g}$ is already adjacent to $v_{4m+2g-1}$.
\end{enumerate}
\begin{figure}[htb]
\vbox{\vskip -.8cm
 \hbox{\centerline{\includegraphics[width=10cm]{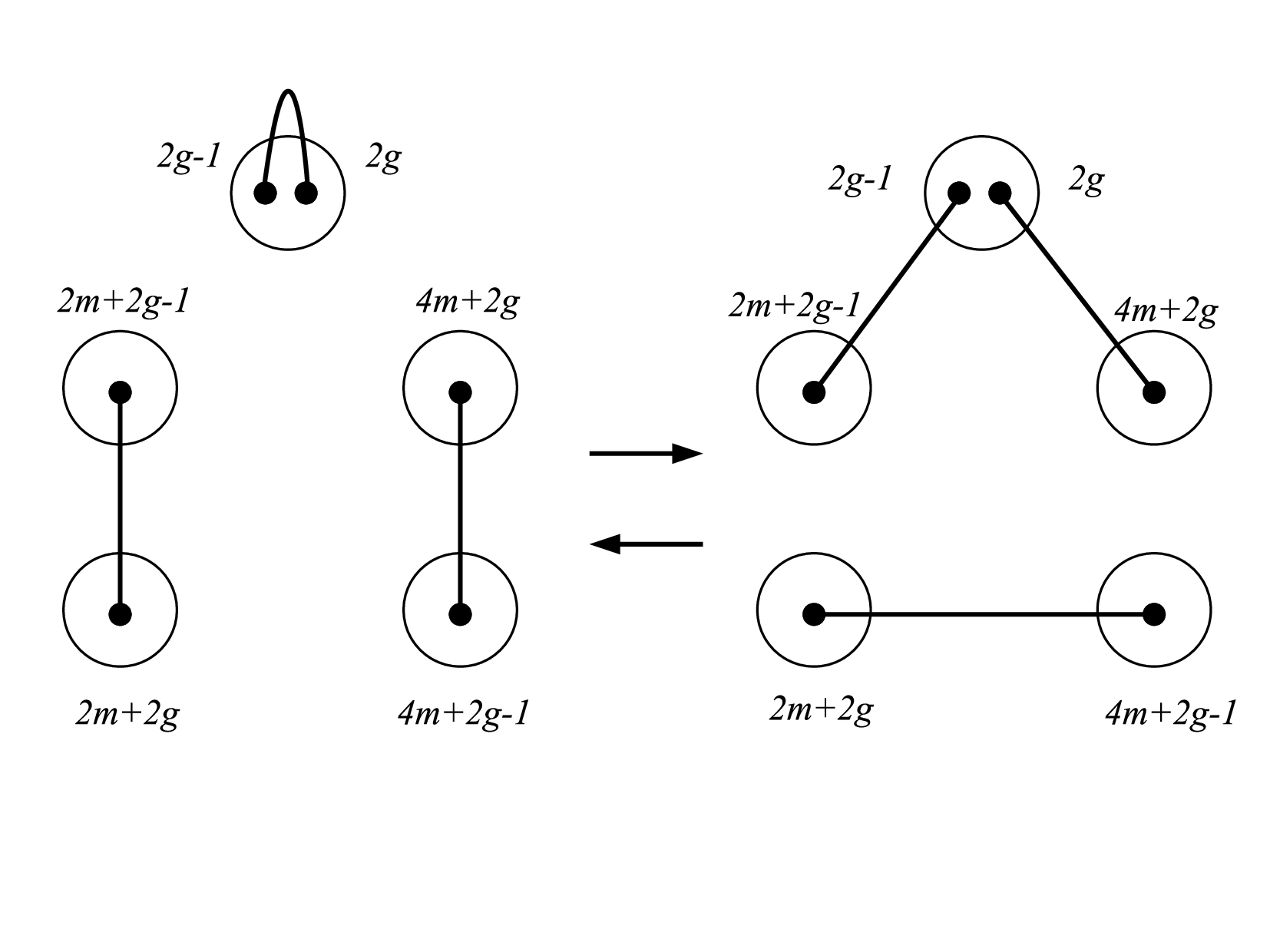}}}

\vskip -2cm \smallskip} \caption{\it  switching for loops}

\lab{f:loop}

\end{figure}

 As in the multiple edge case, the inverse switching has the obvious natural definition. Pick $2m$ points in $i$ and label them $1,\ldots,2m$. Let the point paired to $2g-1$ be labelled $2m+2g-1$ and the point paired to $2g$ be labelled $4m+2g$, for all $1\le g\le m$. Pick another $m$ distinct pairs and label their endpoints $2m+2g$ and $4m+2g-1$ for $1\le g\le m$. The inverse switching replaces $\{2g-1,2m+2g-1\}$, $\{2g,4m+2g\}$ and $\{2m+2g, 4m+2g-1\}$ by $\{2g-1,2g\}$, $\{2m+2g-1,2m+2g\}$ and $\{4m+2g-1,4m+2g\}$. Again, let $v_j$ denote the vertex containing point $j$ for all $2m+1\le j\le 6m$. We say an inverse switching is good  if none of the following conditions holds
 \begin{enumerate}
\item[(d)] the vertices in $\{v_j,\ 2m+1\le j\le 6m\}$ are not all distinct;
\item[(e)] $v_{2m+2g-1}$ is adjacent to $v_{2m+2g}$; or $v_{4m+2g-1}$ is adjacent to $v_{4m+2g}$.
\end{enumerate}

We will use the following auxiliary functions of $\MM\in \M$, some of which have already been defined.

\begin{itemize}
\item $Z_0(\MM)=\sum_{u\le n} m_{u,u}$;
\item $Z_{i,i}(\MM)=Z_0(\MM)-m_{i,i}$;
\item $Z_3(\MM)=\sum_{u\le n} m_{u,u}^2$;
\item $K(\MM)=\sum_{u\le n} (d_u-2)m_{u,u}$;
\item $D(\MM)=\sum_{u\le n} (d_u-2)I_{m_{u,u}\ge 1}$.
\end{itemize}
Let $N_i$ be the expected number of good switchings that can be applied to a random $\P\in \MM_\bit$. There are $m!2^m$ ways to label the endpoints of the $m$ loops at $i$. Potentially there are $M_1^{2m}$ ways to choose the $2m$ pairs, in order, and label their endpoints. Hence, potentially there can be $M_1^{2m}m!2^m$ switchings applied to $\P$. As discussed in the multiple edge case, the cases where the chosen $2m$ pairs are not all distinct contribute a relative error of $O(m^2/M_1)$. Next, we estimate the probability that a random switching is not good when it is applied to a random pairing $\P$.

For (a), for every $2m+1\le j\le 6m$, the probability that $v_j=i$ is $O((d_i-2)/M_1)$. (Note that the subtraction of $2$ is caused by the fact that there are two points in $i$ forming a cycle.) Taking the union bound over $j$, the probability that $v_j=i$ for some $j$ is $O(m(d_i-2)/M_1)$. Next we consider repeated vertices. For every $j=2m+2g-1$, $1\le g\le 2m$, the probability that  $v_j\neq i$ and  $\{j,j+1\}$ forms a loop is $O(Z_{i,i}/M_1)$; hence, the probability that one of the $2m$ pairs forms a loop is $O(mZ_{i,i}/M_1)$. Also, for any two of the $2m$ pairs, the probability that they are both adjacent to a vertex $w$ is $O([d_w]_2/M_1^2)$. Taking the union bound over all $w\le n$ and all choices of two pairs (there are $O(m^2)$ of them) produces the bound $O(m^2 \sum_{w\le n}[d_w]_2/M_1^2)=O(  m^2  { M}_2/M_1^2)$. Hence, the probability that (a) occurs is $O\big((d_i-2)/M_1+mZ_{i,i}/M_1+m^2 { M}_2/M_1^2\big)$.

Let $w$ be a vertex. Similar to the argument in condition (iii) for multiple edges, the probability that $i$ is adjacent to $w$ in a random $\P$ is $O((d_i-2)d_w/M_1)$. (Note that there are at least two points in $i$ that form a loop and and hence unavailable to be paired to a point in $w$.) Conditional on that, then, for any $1\le g\le m$, the probability that either of the points $2m+2g-1$ and $4m+2g$ is in $w$, and (a) does not occur, is $O((d_w-1)/M_1)$. Hence, taking the union bound over $w$ and $g$, the probability that (b) occurs without (a)  is $O(m \sum_{w\le n} [d_w]_2 (d_i-2)/M_1^2)=O(m {M}_2(d_i-2)/M_1^2)$.

By trivial modifications of  the same argument,  the probability that (c) occurs without (a)  is $O(m{M}_2^2/M_1^3)$.

It follows that
\[
N_i=M_1^{2m}m!2^m\Big(1+O(m^2/M_1+(d_i-2)/M_1+mZ_{i,i}/M_1+m^2 { M}_2/M_1^2+m { M}_2(d_i-2)/M_1^2+m{ M}_2^2/M_1^3)\Big).
\]

Next, we estimate $N'_i$, the expected number of good inverse switchings applied to a random $\P\in \C(\MM(0))$. Potentially, there are $[d_i]_{2m}$ ways to choose and label the points $2g-1$, $2g$, $2m+2g-1$, $4m+2g$ for all $1\le g\le m$, and there are $M_1^m$  ways to choose and label the other pairs. Thus, potentially, the number of inverse switchings that can be applied to $\P$ is $[d_i]_{2m}M_1^m$. As before, the proportion of these potential cases where the $m$ randomly chosen pairs $\{2m+2g,4m+2g-1\}$ are not all distinct is $O(m^2/M_1)$. We next estimate the probability that a random switching is not good.

Condition (d) occurs only if some pair $\{2m+2g,4m+2g-1\}$ forms a loop,  or two chosen pairs   use a common  vertex. The probability of the former is $O(mZ_{i,i}/M_1)$. We next  bound the  probability of the latter when a random switching is applied to a random $\P$. There are two subcases. Arguing as before, it is easy to find the probability that  two pairs with the common vertex are both of form $\{2m+2g,4m+2g-1\}$ is $O(m^2{M}_2/M_1^2)$. In the second subcase, one pair uses $i$ and the other pair is of the form $\{2m+2g,4m+2g-1\}$.  We can choose all the points  $2g-1$ and $2g$ in vertex $i$ in advance. The probability that one particular such point  is paired with a point in a given vertex $w\neq i$  is $O(d_w/M_1)$. Conditional upon that, the probability that a given $2m+2g'$ or $4m+2g'-1$ is in $w$ is $O((d_w-1)/M_1)$. Multiplying these by the number $m^2$ of choices for $g$ and $g'$, we see that this subcase contributes the same as the first one.
Hence, the probability for (d) to occur is $O(mZ_{i,i}/M_1+m^2{M}_2/M_1^2)$.

For (e), the analysis is similar to (d), and we easily get $O(m{M}_2^2/M_1^3)$.

We conclude  that
\[
N'_i=[d_i]_{2m}M_1^m(1+O(m^2/M_1+mZ_{i,i}/M_1+m^2{M}_2/M_1^2+m{M}_2^2/M_1^3)).
\]

\subsection{Eliminating  loops} \lab{sec:loops}
Define
$$
\kappa_i(\MM)=m_{i,i}^2/M_1+(d_i-2)/M_1+m_{i,i}Z_{i,i}/M_1+m_{i,i}^2 {M}_2/M_1^2+m_{i,i} { M}_2(d_i-2)/M_1^2+m_{i,i}{ M}_2^2/M_1^3.
$$
Then
$$
\frac{|\C(\MM)|}{|\C(\MM(0))|}=\frac{[d_i]_{2m_{i,i}}}{M_1^{m_{i,i}}m_{i,i}!2^{m_{i,i}}}(1+O(\kappa_i(\MM))),
$$
provided that
\begin{enumerate}
\item[(A4)] $\kappa_i(\MM)\le \xi_1$ for all $i\le n $  such that $m_{i,i}\ge 1$.
\end{enumerate}

For every $\MM$, we can repeatedly switch away all loops in pairings in $\C(\MM_\bit)$ and apply the above estimate for each   ratio as required, and consequently obtain
\be
\frac{|\C(\MM_\bit)|}{|\C(\MMS)|}=
\exp(O(\kappa(\MM)))\prod_{i\le n}\frac{[d_i]_{2m}}{M_1^{m}m!2^m}, \lab{formulaB}
\ee
where $m$ denotes $m_{i,i}$
and
\bea
\kappa(\MM)&=&\sum_{i\le n}\kappa_i(\MM) I_{m_{i,i}\ge 1}\non\\
 &=& Z_3/M_1+D/M_1+Z_3 {M}_2/M_1^2+K {M}_2/M_1^2+Z_0{M}_2^2/M_1^3
 +\sum_{i\le n}m_{i,i}Z_{i,i}/M_1.\lab{kappa}
\eea

\subsection{Combining the switchings to obtain simple pairings}
\lab{s:combining}

Define $\xi(\MM)=\eta(\MM)+\kappa(\MM)$.
 Let $ \M(\xi_1)$ be the set of ${\bf M}\in \M$   for which   $\xi({\bf M})  \le \xi_1$. Obviously every $\MM\in\M(\xi_1)$ satisfies all assumptions (A1)--(A4) . (By considering  the first term in $\eta(\MM)$, we get $ Z_2(\MM)/M_1\le \xi_1$, which implies (A1) and also (A2)  since $Z(\MM)\le Z_2(\MM)$; (A3) and (A4) are implied as $\eta_{i,j}(\MM,m_{i,j})\le \eta(\MM)$ and $\kappa_i(\MM)\le \kappa(\MM)$).
  Thus,  combining~\eqn{formula} and~\eqn{formulaB},
  for every $\MM\in \M(\xi_1)$,
\[
\frac{|\C(\MM)|}{|\C(\MMS)|}=F(\MM)\exp(O(\xi_1)),
\]
where
\be
F(\MM)=\prod_{(i,j)\in \Heavy({\bf M})} \frac{[d_i]_{m_{i,j}}[d_j]_{m_{i,j}} /(m_{i,j}!M_1^{m_{i,j}})}{1+ d_id_j /M_1}\prod_{i\le n}\frac{[d_i]_{2m_{i,i}}}{(2M_1)^{m_{i,i}}m_{i,i}!}
\lab{F}
\ee

From~\eqn{Mtophi} and then~\eqn{formula}, we have
\bea
\frac{|\Phi|}{ \big|\C\big(\MMS\big)\big|}
&=&
\sum_{{\bf M} \in \M }  \frac{ |\C\big({\bf M}  ) |}{ \big|\C\big(\MMS\big)\big|}\non\\
&=&  \frac{ S_1}{ \big|\C\big(\MMS\big)\big|}+\sum_{{\bf M} \in \M( \xi_1 ) }F({\bf M}) \exp\big( O (\xi_1)\big)
\non\\
&=&   \frac{ S_1}{ \big|\C\big(\MMS\big)\big|} + \big(1+ O (\xi_1)\big)S_2\non
\eea
where
$$S_1= \sum_{{\bf M} \in \M\setminus \M(\xi_1) }  |\C\big({\bf M}  ) |,
\qquad S_2=
\sum_{{\bf M} \in  \M(\xi_1) } F({\bf M}).
$$

Hence
\bel{S1S2}
\frac{|\Phi|}{ \big|\C\big(\MMS\big)\big|} (1- S_1/|\Phi| )
= S\big(1+O(1-S_2/S)+   O (\xi_1)\big)
\ee
where
 \bel{S}
 S =  \sum_{{\bf M} \in \M }F({\bf M})= \prod_{1\le i<j<n}  A_{i,j}  \prod_{1\le i\le n}  B_{i}
 \ee
and
\bel{Aij}
A_{i,j}=\sum_{m\ge 0}\frac{[d_i]_m[d_j]_m /(m!M_1^m)}{1+ d_id_j /M_1}, \quad  B_i=\sum_{m\ge 0}\frac{[d_i]_{2m}}{(2M_1)^mm!}.
\ee
Note that the terms $1+ d_id_j /M_1$, for $m=0$ and 1 respectively, in the numerator  of $A_{i,j}$   appear from the case that $m_{i,j}=\bit$, which essentially contributes a factor 1 to the first product in~\eqn{F}.

We will later find  bounds for $   S_1/ |\Phi| $ and $1-S_2/S  $.
 First we analyse $\xi\big({\bf M} \big)$ in order to find a suitable value for $\xi_1$.

\subsection{Bounding $S_1$} \lab{S1}
In this section, our aim is to find a good upper bound, $\xi_2$, on  $S_1/|\Phi|$  for some suitably small value of $\xi_1$. Our final error term will be $O(\xi_1+\xi_2)$. We may view $\xi(\MM)$ as the total of the error bounds for the individual switchings that are relevant to  pairing $\P\in\Phi$ given $\MM(\P)$.     In earlier applications of the switching method to counting graphs with given degrees, the analogue of $\xi_1$ was a bounded multiple of the analogue of $ \ex\,\xi(\MM)$ (viewed in this way).  For those familiar with the argument, a bounded multiple is clearly optimal. This was relatively straightforward in those applications because  the error bound per switching was a simple function of the basic variables being analysed. Roughly speaking, these correspond to $Z$ and $Z_0$. Unfortunately, we cannot afford this luxury in our application because our approach is quite different and we deal with the much more complicated $U$ functions.
Consequently, we content ourselves with $\xi_1$ and $\xi_2$ being approximately the square root of $\ex \, \xi(\MM(\P))$.
That is, our goal is prove that with probability $1-O(\sqrt{\xi})$, $\xi(\MM(\P))\le \sqrt{\xi}$.  We start by evaluating the expectation of each term in $\xi(\MM(\P))$.
Recall the definitions of $U_i$ in~\eqn{tildeUs}.
 We further define
\bean
U_6 &=& \sum_{i\ne j\le n}
\min\{[d_i ]_3[d_j]_3/M_1^2,d_id_j \},\non\\
U_7&=&\sum_{i\ne j\le n}\frac{[d_i]_2}{M_1}\min\{(d_i-2)[d_j]_2/M_1^2,d_j/M_1\}.\non
\eean

\begin{lemma}\lab{lem:Err}
We have the following, where ${\cal H} ={\cal H}(\MM(\P))$  is   defined  above~\eqn{formula}  and all functions are of a  pairing $\P$ taken u.a.r.\ from $\Phi$. Here $m_{i,j}$ refers to the entry of $\MM(\P)$.
\bean
&&\ex Z=O({U}_2);\quad \ex Z_2=O({M}_2^2/M_1^2);\quad \ex \sum_{\cal H} m_{i,j}Z_{i,j} =O({U}_2^2);\\
&&\ex \sum_{\cal H} m_{i,j}(W_{i,j}+W_{j,i}) =O({ U}_3);\quad \ex \sum_{\cal H} m_{i,j}(d_i-2) =O({U}_4);\\
&&\ex  \sum_{\cal H} m_{i,j} (R_{i,j}/d_i + R_{j,i}/d_j)=O({U}_5);\quad \ex \sum_{\cal H} m_{i,j}Q_{i,j} =O({ U}_2{ U}_6);\\
 &&\ex \sum_{\cal H} m_{i,j}   (m_{i,i}/d_i+m_{j,j}/d_j) = O(U_7) ;\quad \ex ZZ_0=O(U_2{M_2}/M_1);\\
&&\ex Z_0=O({M_2}/M_1);\quad \ex Z_3=O({M}_2/M_1+{M}_4/M_1^2);\\
&&\ex K=O({M}_3/M_1);\quad \ex D  =  O(U_1);\quad  \ex \sum_{1\le i\le n}m_{i,i}Z_{i,i}=O(M_2^2/M_1^2).
\eean

\end{lemma}

\proof
  An upper bound on   $\ex Z(\P)$ is obtained as follows. First, note that if $Y_{u.v}$ is the multiplicity of the edge $uv $ in $\P$,
 $$
  \ex( Y_{u,v}I_{Y_{u,v}\ge 2})\le \ex \big(\min\{ [Y_{u,v}]_2,  Y_{u,v} \}\big) \le \min\{\ex  [Y_{u,v}]_2,\ex Y_{u,v} \}.
$$
The number of locations for two non-loop pairs in parallel is at most $\sum [d_u]_2[d_v]_2$,
summed over all   vertices $u<v$, and similarly $\sum d_ud_v$ for just one pair.
 Since all of the $M_1$ points are uniformly at random (u.a.r.) paired in $\P$, for any constant integer $k>0$, the probability of a given set of $k$ pairs occurring in $\P$ is
$$
\frac{\prod_{i=k}^{M_1/2-1}(M_1-2i-1)}{\prod_{i=0}^{M_1/2-1}(M_1-2i-1)}=O(M_1^{-k}).
$$
  Hence, recalling the definition of $U_2$
from~\eqn{tildeUs},
we have
\bel{EZP}
\ex Z(\P) =  \sum_{1\le u<v\le n} \ex( Y_{u,v}I_{Y_{u,v}\ge 2})=O(U_2).
\ee
Similarly,
\bel{EZ2P}
\ex Z_2(\P) =   \ex\sum_{1\le u<v\le n }Y_{u,v}^2 I_{Y_{u,v}\ge 2}\le 2\sum_{1\le u<v\le n } \ex( Y_{u,v}( Y_{u,v}-1 ))=O(M_2^2/M_1^2).
\ee
 We apply  a similar argument to the other terms in the lemma.
Firstly,
$$
\ex \sum_{\cal H}  m_{i,j} {Z_{i,j}} = O(U_2^2).
$$
 To see this, note that for any pair of   vertices $(u,v)$ of concern other than $(i,j)$, the bound
$\min\{d_ud_v/M_1,d_u^2d_v^2/M_1^2\}$ is still valid for   $\ex  Y_{u,v}I_{Y_{u,v}\ge 2} $ even when a given value of $m_{i,j}$ is conditioned upon. Hence the expression is   bounded above by $(\ex Z(\P))^2$.

Next, we show
\be
\sum_{\cal H} m_{i,j}(W_{i,j}+W_{j,i})=O(U_3).\lab{U3}
\ee
First notice that
\bea
\ex\sum_{\cal H} Y_{i,j}W_{i,j}&=&\ex\Big(\sum_{1\le i<j\le n} m_{i,j} W_{i,j} I_{Y_{i,j}\ge 2}\Big)\non\\
&=&\ex\Big(\sum_{1\le i<j\le n} Y_{i,j}I_{Y_{i,j}\ge 2}\sum_{w\notin\{i,j\}} (d_w-2) I_{Y_{i,w}\ge 2}\Big).\lab{productEx}
\eea
Given any value of $Y_{i,j}\ge 2$, the conditional expectation of
$\sum_{w\notin\{i,j\}} (d_w-2) I_{Y_{i,w}\ge 2}$ is always bounded by
$\sum_{w\le n} (d_w-2) O\big(\min\{[d_i-2]_2[d_w]_2/M_1^2, 1\}\big)$
using the fact that
\bel{PrMult}
\pr( Y_{i,w}\ge 2 \mid Y_{i,j}\ge 2 ) = O\big(\min\{[d_i-2]_2[d_w]_2/M_1^2, 1\}\big).
\ee
Hence, we can separate the product inside the expectation in~\eqn{productEx} and bound its expectation asymptotically  (within a constant factor) by the product of the two expectations. We have already shown that
$$
\ex \Big(\sum_{1\le i<j\le n} Y_{i,j}I_{Y_{i,j}\ge 2}\Big)=\ex Z =O(U_2).
$$
 Recalling the definition of $U_3$ from~\eqn{tildeUs}, we have $\ex\sum_{\cal H} Y_{i,j}W_{i,j}=O(U_3)$. By swapping the labels of $i$ and $j$ and noting that in the definition of $U_3$, $i$ and $j$ are not ordered, we also have $\ex\sum_{\cal H} Y_{i,j}W_{j,i}=O(U_3)$, and hence~\eqn{U3}.

It is straightforward to bound the expectations of all the other terms in the lemma  in a similar fashion. \qed\ss

By Lemma~\ref{lem:Err} and the definition  of $\xi(\MM)=\eta(\MM)+\kappa(\MM)$ where $\eta(\MM)$ and $\kappa(\MM)$ are given in~\eqn{eta} and~\eqn{kappa}, we have
$\ex\,\xi(\P) = O(\xi_0)$ where
\bea
\xi_0 &=& U_5+ U_7+
\frac{   U_2^2+ U_3+  U_4+   U_1}{M_1 }
+\frac{  U_2  M_2+   U_4   M_2+   U_2   U_6}{M_1^2}
+\frac{     U_2    M_2^2}{M_1^3 }
\non \\
 &&+\frac{  M_2}{M_1^2}+\frac{  M_4+  M_3 M_2+  M_2^2}{M_1^3}+\frac{  M_4  M_2+  M_2^3}{M_1^4}.\non
\eea
Note that, by elementary considerations,
 \bean
 &&  U_7=O( M_3 M_2/M_1^3); \quad \frac{ U_4}{M_1}=O(  M_3  M_2/M_1^3); \quad \frac{ U_2 M_2}{M_1^2}=O( M_2^3/M_1^4);\\
 &&   M_4\le   M_3  M_2;\quad  M_4  M_2\le  M_2^3;\quad
    M_2^2/  M_1^3 \le  M_2 /  M_1^2 +  M_2^3/M_1^4.
 \eean
Moreover, by the hypothesis of Theorem~\ref{thm:general}  that
$\xi=o(1)$, we have (taking the first option in the min functions) $ U_2  U_6/M_1^2=O(  M_2^2  M_3^2/M_1^6)=o( M_2  M_3 /M_1^3)$.
Thus
$\xi_0 = O(\xi)$ and we have  $\ex\xi(\P)=O(\xi)$.

 Now we set $\xi_1$ in the previous sections to be $\sqrt{\xi}$, which is $o(1)$.
This definition determines the precise set $\M( \xi_1)$ which we have been dealing with since Section~\ref{s:combining}, and ensures that $\xi_1 = o(1)$ as required by (A1--A4).

Recalling the definition of $S_1$ above~\eqn{S1S2}, the  following comes immediately from Lemma~\ref{lem:Err} using Markov's inequality.

\begin{cor} \lab{cor2:xiMax}  Assume that  the  hypothesis of Theorem~\ref{thm:general} holds. With probability $1-O\big(\sqrt{\xi}\big)$, ${\bf M}(\P)\in \M\big(\sqrt{\xi}\big)$.
\end{cor}
It follows from this corollary that
\bel{S1phi}
\frac{ S_1}{|\Phi|} =
 \pr\big({\bf M}(\P)\notin \M(\sqrt{\xi})\big)= O\big( \sqrt{\xi}\big).
\ee
\subsection{Bounding $S-S_2$ }
Next, as might be foreseen from~\eqn{S1S2}, we wish to bound $1-S_2/S  $. Define a  probability space $\Omega^*$ by equipping $\M$ with a new probability function, in which $\pr({\bf M})$ is proportional to the ``weight'' $F({\bf M})$ defined in~\eqn{F}. Then, noting that the total weight is $S$, we have  $1-S_2/S =\pr\big({\bf M}\notin \M(\sqrt{\xi})\big)= \pr(\xi({\bf M})>\sqrt{\xi})$ by definition of $\M(\sqrt{\xi})$. (Recall that we set $\xi_1=\sqrt \xi $ in Section~\ref{S1}.)   Next, observe that $\Omega^*$ is a product space with each $m_{i,j}$ chosen independently at random from the distribution of a random variable $X_{i,j}$ defined as follows, where the normalising factors $A_{i,j}$ and $B_i$  are given  in~\eqn{Aij}.   Let $\pr(X_{i,j}=\bit) =   A_{i,j}^{-1}$ for $i<j$,     and
\bean
\pr(X_{i,j}=m) &=&  A_{i,j}^{-1}\frac{[d_i]_m[d_j]_m /(m!M_1^m)}{1+ d_id_j /M_1} \quad (i<j,\ m\ge 2)\\
\pr(X_{i,i}=m)&=& B_i^{-1}\frac{[d_i]_{2m}}{m!(2M_1)^m} \quad (m\ge 0).
\eean
Clearly $\pr\big( X_{i,j}\ne \bit\big)=O\big([d_i]_2[d_j]_2/M_1^2\big)$. Let   $ \la_{i,j}  =d_id_j/M_1$. Since $\pr(X_{i,j}=m)/\pr(X_{i,j}=m-1)\le \la_{i,j}/m$ when $m\ge 3$, and this is the corresponding ratio for the Poisson variable $\po(\la_{i,j})$, it follows that $X_{i,j}$ in $\Omega^*$ is stochastically dominated by $\po(\la_{i,j})$ (recalling that $\bit$ is treated as 0 in numerical functions). Hence     $\ex  X_{i,j} = O(d_id_j/M_1)$.
  Similarly, $X_{i,i}$ is stochastically dominated by $\po(\mu_i)$ where $\mu_i =[d_i]_2/2M_1$.

We next show that the expected value of $\xi({\bf M})$ in $\Omega^*$ is $O(\xi)$. The general idea is to show that for each auxiliary function $f\in\{Z,Z_2,\ldots, K, D\}$ defined in Sections~\ref{sec:auxiliary} and~\ref{sec:loops}, the expected value $\ex f({\bf M})$ for ${\bf M}\in \Omega^*$ is close to that in Lemma~\ref{lem:Err} for $\MM(\P)$ where $\P$ is a random pairing in $\Phi$.    We first verify  this in detail  for $f=Z$ and $f=Z_2$.

By definition, $\ex Z\le \sum_{i<j}\ex X_{i,j}=\sum_{i<j}O(d_id_j/M_1)$. Moreover, since $X_{i,j}$ is never equal to 1,
  $\ex  X_{i,j}\le \ex[ X_{i,j}]_2$. Using the domination by Poisson, this is $  O([d_i]_2[d_j]_2/M_1^2)$.  Thus $\ex Z\le \sum_{i<j}\ex[ X_{i,j}]_2 =\sum_{i<j} O([d_i]_2[d_j]_2/M_1^2)$.
   Similarly, by the definition of $Z_2$, we have $\ex Z_2\le \ex Z+ \sum_{i<j} \ex [X_{i,j}]_2=\sum_{i<j} O([d_i]_2[d_j]_2/M_1^2)$. Hence,  recalling the definition of $U_2$ from~\eqn{tildeUs},  we have
$$
\ex Z(\MM) = O(U_2),\quad \ex Z_2(\MM)=  O(M_2^2/M_1^2).
$$
 A similar argument applies to the other terms in $\xi$. For instance, we may bound $\ex \sum  m_{i,j}  W_{i,j}$ by the summation of $[d_w]_3[d_i]_2^2[d_j]_2$ over ordered triples $(w,i,j)$. To obtain the same error term as before, note that $[d_i]_2^2=O([d_i]_4+d_i)$, and we may obtain the other terms in the $\min$ functions in $U_3$ using arguments analogous to those used for the case of random $\P$. The remaining details required for showing $\ex \xi({\bf M})=O(\xi)$ are straightforward.  In particular, note that $m_{i,j}$ and $Z_{i,j}$ are, by design, independent, which makes it easy to write the expected value of $m_{i,j}Z_{i,j}$. (This is why we use $Z_{i,j}$ rather than $Z$.) It then follows, by Markov's inequality, that in $\Omega^*$, $\pr(\xi({\bf M})>\sqrt{\xi})=O(\xi_0/\sqrt{\xi})=O(\sqrt{\xi})$, and thus by the same argument as before, $1-S_2/S=\pr\big({\bf M}\notin \M(\sqrt{\xi})\big)=O(\sqrt{\xi})$.
 Combining  this with~\eqn{S1phi} in~\eqn{S1S2} produces
\bel{donedoubles}
  \frac{|\Phi|}{ \big|\C\big(\MMS\big)\big|} = S\big(1+O(\sqrt{\xi})\big).
\ee

\subsection{Estimating $S$}

Here we obtain a much more user-friendly version of the function $S$ from~\eqn{S}. Note that the extra error term $ M_1^{-1}$ makes no difference if  $\Delta\ge 3$ since then $U_1>0$ and $\xi\ge 1/M_1$. If $\Delta\le 2$ then we could slightly modify the following lemma to go further, but in this case the enumeration problem is anyway easily solved by other means.
 \begin{lemma}\lab{lem:S} Assume  $\xi=o(1)$  and let $\la_{i,j}=d_id_j/M_1$. Then

 \[
 S=(1+O(\xi+M_1^{-1}))\exp\left(\frac{M_1}{2}-\frac{M_2}{2M_1}+\frac{M_3}{3M_1^2}-\frac{3}{4}-\sum_{1\le i<j\le n} \log(1+\la_{i,j}) \right).
 \]
 \end{lemma}

\proof
 We start by analysing $A_{i,j}$. Recall that
\bel{Aijnew}
A_{i,j}(1+ d_id_j /M_1 )=\sum_{m\ge 0}\frac{[d_i]_m[d_j]_m}{m!M_1^m}.
\ee

   Let $\la_{i,j}=  d_id_j/M_1$.  We only need to  consider   the terms with $m\le m_0=\max\{\log^2 M_1,C\la_{i,j}\}$ for some sufficiently large constant $C>0$.  This is because elementary arguments, for instance considering ratios of successive terms, easily show that the terms with $m>m_0$ contribute a relative proportion $O(M_1^{-C'})$ of the total summation, where $C'\to\infty$ as $C\to\infty$.     By Lemma~\ref{lem:tau},   the maximum degree $\Delta$ is $o(M_1^{3/5})$.

 We will at first obtain two different formulae, depending on the size of $d_id_j$. We are able to put the second formula into a form that is valid in both cases for an appropriate choice of the split between cases. First we need some observations about the equation
\bel{stirling1}
 [d]_m/d^m=\sum_{k =0 }^{ m-1} s(m,m-k) d^{-k}
\ee
where $s(m,m-k)$ is a Stirling number of the first kind. By definition,
$s(m,m)=1$ and
$$
s(m,m-k)=\sum_{1\le b_1<b_2<\cdots<b_k<m}\prod_{i=1}^k (-b_i)  =\sum_{1\le b_k< m}  (-b_k)s\big(b_k,b_k-(k-1)\big)\quad \ \mbox{for all}\ 1\le k<m.
$$
Hence, by induction on $k$, $s(m,m-k)$ is a polynomial $P_k$ in $m$ of degree at most $2k$, defined by
$$
P_k(m)=\sum_{b=1}^{m-1} -bP_{k-1}(b),\qquad P_0=1
$$
using the standard   formula for $\sum_{b=1}^{m-1} b^r$.
Note that this is valid even for $k\ge m$, when it evaluates to 0.

The first part of the above equation also gives
$ |s(m,m-k)|\le \Big(\sum_{i=1}^{m}  i\Big)^k<m^{2k}$. Hence
for any $d$ and $m$ with $m^2/d\le 1/2$  say,  and for any integer $0<t\le m$,

\be
 \sum_{k= t}^{m-1} |s(m,m-k)| d^{-k}=O(m^{2t}/d^t). \lab{tailBound}
\ee

\smallskip

\noindent{\em Case 1:\/}  $d_id_j \ge M_1^{5/6}$.
\smallskip

Our main object is to analyse $[d_i]_m[d_j]_m$ in~\eqn{Aijnew}, for which we will use~\eqn{stirling1}.
Since  $d_j\le \Delta=o(M_1^{3/5})$,   we have $d_i =\Omega(M_1^{5/6}/d_j)=\Omega(M_1^{7/30})$.    Hence, if $m\le \log^2 M_1$, then  $m^2/d_i=o(1/M_1^{1/5} )$,    and similarly for $d_j$.  On the other hand, if $   \log^2 M_1<m\le m_0$ then  $m=O( \la_{i,j})$  and hence  $m^{2  }/d_i = O(d_i d_j^{2 }/M_1^{2 })= O(\Delta^{3   }/M_1^{2 })=  O(M_1^{-1/5})$.
In both cases, we have  by~\eqn{stirling1} and~\eqn{tailBound} that for fixed $u$
 \bel{stirling2}
[d ]_m/d ^m = \sum_{r=0}^{u-1}s(m,m-r)d^{-r}  +O(m^{2u }/d^{u })= \sum_{r=0}^{u-1}s(m,m-r)d^{-r}  +�O(M_1^{-u/5}).
\ee
Similarly, since $s(m,m-k)$ is a polynomial $P_k$ in $m$ of degree at most $2k$,
\bel{neweq}
s(m,m-r)d_i^{-r} s(m,m-w)d_i^{-w} \le m^{2r+2w}/(d_i^rd_j^{w})
\ee
   for every  fixed $r$ and $w$. Thus for $c_1=5K-1$ ($K$ fixed) we have
$$
[d_i]_m[d_j]_m  =  O(M_1^{-K})+ d_i^md_j^m \sum_{r=0}^{c_1}\sum_{w=0}^{c_1-r} s(m,m-r)d_i^{-r}s(m,m-w)d_j^{-w}.
$$
   Rewriting the polynomials $P_k$ in  terms of the falling factorials $[m]_t$  using Stirling numbers of the Second kind,  we have
for each $r$ and $w$ with $r+w\le c_1$, that $s(m,m-r)s(m,m-w)=\sum_{t=0}^{2c_1}a_{r,w,t}[m]_t$ for some absolute constants $a_{r,w,t}$ where $a_{0,0,0}=1$. Hence
\bel{didj}
[d_i]_m[d_j]_m  =    O(M_1^{-K})+d_i^md_j^m \sum_{r=0}^{c_1}\sum_{w=0}^{c_1-r}\sum_{t=0}^{2c_1} a_{r,w,t}[m]_td_i^{-r}d_j^{-w}.
\ee

We may now rewrite~\eqn{Aijnew}, recalling that terms with $m>  \max\{\log^2 n, C\la_{i,j}\}$ can be ignored in~\eqn{Aijnew}. Noting that the new terms introduced in the following are similarly negligible, we have
  \bea
A_{i,j}(1+ \la_{i,j})&=& \sum_{m\ge 0}\frac{    \la_{i,j}^m }{m! }\bigg( O(M_1^{-K})+\sum_{r=0}^{c_1}\sum_{s=0}^{c_1-r}\sum_{t=0}^{2c_1} a_{r,w,t}[m]_td_i^{-r}d_j^{-w}\bigg). \non\eea
Thus, using
\bel{expsum}
\sum_{m\ge 0}  [m]_t\la_{i,j}^m/m! =\la_{i,j}^t \exp(\la_{i,j}) = d_i^td_j^tM_1^{-t}\exp(\la_{i,j})
\ee
 we get
\bel{UPPER0} A_{i,j}(1+ \la_{i,j})\exp(-\la_{i,j}) =  O(M_1^{-K}) + \sum_{r=0}^{c_1}\sum_{w=0}^{c_1-r}\sum_{t=0}^{2c_1}  a_{r,w,t}  d_i^{t-r}d_j^{t-w}/M_1^t
\ee
and hence
\bel{UPPER} \log\big(A_{i,j}(1+ \la_{i,j})\big) =  \la_{i,j} +  O(M_1^{-K}) + \log\sum_{r=0}^{c_1}\sum_{w=0}^{c_1-r}\sum_{t=0}^{2c_1}  a_{r,w,t}  d_i^{t-r}d_j^{t-w}/M_1^t.
\ee

\smallskip

\noindent{\em Case 2:\/}  $d_id_j\le  M_1^{1-\eps}$ for some $0<\eps\le1/6$.
\smallskip

In this case, $[d_i]_m[d_j]_m\le  M_1^{m-m\eps}$. Hence, in the summation~\eqn{Aij}  defining $A_{i,j}$,  the sum of terms for $m> c_2$ for any constant $c_2\ge \lceil K/\eps\rceil -1$ is  $O(1/M_1^{K})$ (with $K$ as in Case~1).  That is,
$$
A_{i,j}(1+\la_{i,j})=O(M_1^{-K})+\sum_{m= 0}^{c_2} [d_i]_m[d_j]_m /(m!M_1^m),
$$
and thus it is straightforward to verify that
\bel{LOWER}
\log\big(A_{i,j}(1+ \la_{i,j})\big) = O(M_1^{-K})+\phi_{c_2}\left(\log\sum_{m= 0}^{c_2} [d_i]_m[d_j]_m /(m!M_1^m)\right)
\ee
where $\phi_{c_2}$ truncates the expansion of the logarithm of the summation, deleting any terms containing   $M_1^{-u}$ for $u>c_2$.

We next rewrite~\eqn{LOWER} into a form that we show is equivalent to~\eqn{UPPER} when $c_2$ is large enough. (This equivalence could alternatively be shown by direct but tedious---especially in the case of~\eqn{UPPER}---computation for any particular value of $K$.)
It is a quite subtle aspect of our argument that the terms truncated by $\phi_{c_2}$ must not be included! (They make no difference for Case~2 but would spoil Case~1.)
Fix any positive constant $c_2$. Recalling that $s(m,m-k)$ is a  polynomial $P_k(m)$ of degree at most $2k$,  we may start with~\eqn{stirling1} and apply the argument leading to~\eqn{didj}
but retaining all the terms in the expansion~\eqn{stirling2} for $ r \le  c_2 -1  $.  Recalling that the polynomials  $s(m,m-r)$    do
their job even for  $m\le r$, and  in~\eqn{LOWER} we only consider $m\le c_2$ (implying the error in~\eqn{stirling2} in this case becomes zero),   we have for $m\le c_2$
$$ [d_i]_m[d_j]_m = \sum_{t=0}^{2c_2-2} (d_id_j)^mQ_t(d_i^{-1},d_j^{-1}) [m]_t$$
where, with $a_{r,w,t}$ as in~\eqn{UPPER},
$$Q_t(d_i^{-1},d_j^{-1})=\sum_{r=0}^{c_2-1}\sum_{w=0}^{c_2-1}a_{r,w,t} d_i^{-r}d_j^{-w}.$$

  Now let us consider, as a formal power series in $z$,
\bean
\log\sum_{m= 0}^{c_2} \frac{[d_i]_m[d_j]_m }{m!}z^m
&=& \log\sum_{m= 0}^{c_2} \sum_{t=0}^{m} \frac{(d_id_j)^m}{m!}Q_t(d_i^{-1},d_j^{-1}) [m]_t z^m\\
&=& \log\left(\sum_{t=0}^{c_2} Q_t(d_i^{-1},d_j^{-1}) \sum_{m= t}^{c_2}  \frac{(d_id_j)^m}{m!} [m]_t z^m\right)\\
&=& \log\left(\sum_{t=0}^{c_2} Q_t(d_i^{-1},d_j^{-1}) \Big(\exp(zd_id_j) (zd_id_j)^t +O(z^{c_2+1})\Big)\right).
\eean
 where, in the first step, we note that $[m]_t=0$ if $t>m$, and in the last step,  $O()$ is used in the formal power series sense and follows from the formula for $\sum_{m\ge 0}\frac{x^m}{m!} [m]_t.$
The latter expression is
$$
\log\left(\exp(zd_id_j)\sum_{t=0}^{c_2} Q_t(d_i^{-1},d_j^{-1})(zd_id_j)^t )  \right)    +O(z^{c_2+1})$$
$$
=zd_id_j + \log\left( \sum_{t=0}^{c_2} Q_t(d_i^{-1},d_j^{-1})(zd_id_j)^t )  \right)  +O(z^{c_2+1}).
$$
Hence, substituting $M_1^{-1}$ for $z$,
\bel{LOWER2}
\phi_{c_2}\left(\log\sum_{m= 0}^{c_2} [d_i]_m[d_j]_m /(m!M_1^m)\right)
= d_id_j/M_1 + \phi_{c_2}\left(\log \sum_{t=0}^{c_2} Q_t(d_i^{-1},d_j^{-1})( d_id_j/M_1)^t )  \right).
\ee
Note that the right hand side has the same terms as~\eqn{UPPER} from Case~1, but with a different cut-off. We are free to choose   $c_2\ge 2c_1=10K-2$ (provided $c_2\ge \lceil K/\eps\rceil-1$), in which case every term in~\eqn{UPPER0} appears in the right hand side of~\eqn{LOWER2} inside the logarithm. As we noted using~\eqn{tailBound} in Case~1, including a bounded number of extra terms  in the expansion   in Case~1 adds an error of order $O(M_1^{-K})$.  Since $a_{0,0,0}=1$, such extra terms would contribute $O(M_1^{-K})$ in~\eqn{UPPER}. Assuming that $\eps=1/6$ and $K\ge 1$, this shows that~\eqn{LOWER} is also valid in Case~1 for   $c_2=   10K-2$   (though some terms in the formula will be dominated by the error term).

We can approximate the simpler function $B_i$  in a similar way. Recall that
$$
B_i=\sum_{m\ge 0}\frac{[d_i]_{2m}}{(2M_1)^mm!}.
$$
When $d_i> M_1^{4/5}$, we may terminate the expansion of $[d_i]_{2m}/d_i^{2m}$  at $m=c_3$ where  $c_3=5K -1$  by incorporating an $O(M_1^{-K})$ error term, and this leads to
\be
\log B_{i} = \frac{d_i^2}{2M_1} +  O(M_1^{-K}) + \log\sum_{s=0}^{c_3}\sum_{t=0}^{2c_3}  b_{s,t} d_i^{t-s}/M_1^t,\lab{UPPER1}
\ee
analogous to~\eqn{UPPER}, where $b_{s,t}$ are absolute constants independent of $d_i$ with $b_{0,0}=1$. On the other hand, for $d_i\le M_1^{4/5}$, we may terminate the summation in $B_i$ at $m=c_4$ such that $c_4\ge 5K -1$   by incorporating an $O(M_1^{-K})$ error, and this yields
\be
\log B_{i} = O(M_1^{-K})+\phi_{c_4}\left(\log\sum_{m= 0}^{c_4} [d_i]_{2m}/(m!(2M_1)^m)\right).\lab{LOWER1}
\ee
With the same argument as for $A_{i,j}$,  with any choice of fixed $c_4\ge  2 c_3=10K-2$, all significant terms in~\eqn{UPPER1} appear  in~\eqn{LOWER1} and all terms in~\eqn{LOWER1} not appearing in~\eqn{UPPER1} are insignificant when $d_i> M_1^{4/5}$.  Thus~\eqn{LOWER1} is valid for any $d_i$ satisfying the conditions of the lemma.

We now choose $K=3$ in both~\eqn{LOWER} and~\eqn{LOWER1} and consequently, to ensure the equivalence shown above,  $c_2=c_4=\ctwo$,   to obtain from~\eqn{S}
 \bean
 S&=&
 \exp\left(\sum_{1\le i<j\le n} \log A_{i,j} + \sum_{1\le i\le n} \log B_i \right)\\
 &=&\exp\left(\psi({\bf d})-\sum_{1\le i<j\le n} \log(1+\la_{i,j})+O(n^2M_1^{-3})\right),
 \eean
 where
 \[
 \psi({\bf d})=\sum_{1\le i<j\le n}\phi_{\ctwo}\left(\log\sum_{m= 0}^{\ctwo} [d_i]_m[d_j]_m /(m!M_1^m)\right)+\sum_{1\le i\le n}\phi_{ \ctwo }\left(\log\sum_{m= 0}^{ \ctwo } [d_i]_{2m}/(m!(2M_1)^m)\right).
 \]
  Noting that $M_1=\Omega(n)$ by our assumption that $d_1\ge 1$, the error term $n^2M_1^{-3}$ is $O(M_1^{-1})$. It only remains to show that
  \be
  \psi({\bf d})=\frac{M_1}{2}-\frac{M_2}{2M_1}+\frac{M_3}{3M_1^2}-\frac{3}{4}+O(\xi+M_1^{-1}).\lab{final}
  \ee
The functions $\phi_{\ctwo}$  is simply a truncation  of the expansion of the logarithm. The result, using Maple for example, is~\eqn{final}, after neglecting all terms that are dominated by $\xi$ or $M_1^{-1}$. This completes the proof of Lemma~\ref{lem:S}.\qed \ss

 We observe that, somewhat surprisingly,  in expanding~\eqn{LOWER}, all terms of the form $d_i^rd_j^w/M_1^t$  with $r+w>t+1$ apparently  disappear, and similarly all terms in~\eqn{LOWER1}  of the form $d_i^w/M_1^t$ with $w>t+1$ disappear.
 The question of finding a proof of this claim  was posed, in the case of~\eqn{LOWER},  in a problem session at a meeting in Oberwolfach\footnote
 {``Enumerative Combinatorics''\!\!, Oberwolfach Workshop ID 1410, 2--8 March, 2014. Organisers  M. Bousquet-M{\'e}lou, M. Drmota, C. Krattenthaler and M. Noy} and immediately solved independently by each of I. Gessel, G. Schaeffer and R. Stanley. Gessel did the same for~\eqn{LOWER1}.
Using these facts, one can reduce the amount of computation required, by first showing that the terms with $4\le m\le 28$ from~\eqn{LOWER} and~\eqn{LOWER1} are dominated by $\xi+O(M_1^{-1})$. Here, it helps to observe that $\sum_{r<w} d_i^rd_j^w \le (M_2^*)^{(r+w)/2}$ for fixed $r,w\ge2$.

\subsection{Completing the proof}

Note that our goal is to prove that
\be
\pr(\G^*(n,{\bf d})\ \mbox{is simple})=(1+O(\sqrt{\xi}+M_1^{-1}))\exp\left(-\frac{M_1}{2}+\frac{M_2}{2M_1}-\frac{M_3}{3M_1^2}+\frac{3}{4}+\sum_{i<j}\Big(\log(1+d_id_j/M_1)\Big)\right)\lab{finalprobability}
\ee
and then our theorem follows by~\eqn{equiv}.

Since $\pr(\G^*(n,{\bf d})\ \mbox{is simple})$ equals $|\C\big(\MMS\big)|/|\Phi|$ and by~\eqn{donedoubles}, it equals $(1+O(\sqrt{\xi}))S^{-1}$. Thus~\eqn{finalprobability} follows by Lemma~\ref{lem:S} (and noting that $\xi=O(\sqrt{\xi})$). This completes the proof of Theorem~\ref{thm:general}.

\section{Bounds on the error}
\lab{sec:boundxi}

 Let $\xi$ be defined as in Theorem~\ref{thm:general}.
Due to the complexity of its definition via the $U_k$, we present in this section several simpler bounds on $\xi$, which are tight in many situations. We first prove Lemma~\ref{lem:xisimple1} presented in Section~\ref{sec:results}.\ss

\no {\bf Proof of Lemma~\ref{lem:xisimple1}.\ }
We can take the first item in each minimum function in the definition of $U_k$. So,
$U_2\le M_2^2/M_1^2$,
$U_3\le M_2 M_3 M_4 /M_1^4$,
$U_4\le M_2  M_3 /M_1^2$,
$U_5\le M_2^2 M_3  /M_1^4$,
$U_1\le M_3  /M_1$.
This  immediately gives the claimed bound on $\xi$ but with an extra term $M_2 M_3 /M_1^3$. However, this term is
  dominated by $M_3/M_1^2+M_2^2M_3/M_1^4$. This completes the proof for part (a).

If we have $\Delta=O(M_1^{1/2})$, then $M_2 M_3 M_4/M_1^5$  is
dominated by $M_2^2 M_3/M_1^4$, since $ M_4=O(d_1^2 M_2)=O(M_1 M_2)$. This bound on $\xi$ is tight within a constant factor, because for such $\Delta$  the first item in the minimum function in each    $U_k$ dominates the second.
 \qed

 The following corollary of Lemma~\ref{lem:xisimple1}  is intended for use when there are not too many vertices with degree less than $3$.
\begin{cor}\lab{cor2:simplexi} Putting $M^*_i=M_i+M_1$ for $i=2,3$, we have
\begin{enumerate}
\item[(a)] $\xi=O\big(M_3^* (M_2^*)^2/M_1^4+M_4M_3M_2/M_1^5\big)$.
\item[(b)] If  $\Delta=O(\sqrt{M_1})$  then $\xi=O\left(M_3^* (M_2^*)^2/M_1^4\right)$.
\end{enumerate}
\end{cor}

\proof It is easy to see that $(M_2+M_3)/M_1^2$ is bounded by $M_3^*(M_2^*)^2/M_1^4$. It is also easy to see that $M_2\le M_3^*$ which eliminates $M_2^3/M_1^4$ from the bound in Lemma~\ref{lem:xisimple1}. Using the Cauchy inequality, we have $M_2^2=O(M_3^*M_1)$ which eliminates $M_2^4/M_1^5$.  Now part (a) immediately follows from Lemma~\ref{lem:xisimple1}(a) and part (b) follows from Lemma~\ref{lem:xisimple1}(b).\qed \ss

The next result follows immediately.
\begin{cor}\lab{cor3:simplexi}
Suppose $\Delta=O(\sqrt{M_1})$, $M_1=O(M_2)$ and $M_1=O(M_3)$. Then
$
\xi=\Theta\left(M_3M_2^2/M_1^4\right).
$
\end{cor}
\no{\bf Remark}.  This bound    on $\xi$   is tight within a constant factor, since by Lemma~\ref{lem:xisimple1}(b), we have $\xi=\Omega(M_3M_2^2/M_1)$, and since $M_1=O(M_2)$ and $M_1=O(M_3)$ imply   $M_2^*=\Theta(M_2)$ and $M_3^*=\Theta(M_3)$.
\smallskip

We now present some results more useful when $\Delta$ is large.  Without loss of generality we may assume that $d_1\ge  \cdots \ \ge  d_n\ge 1$.  First, choose $1\le h < n$ and define $H_k=\sum_{i\le h}[d_i]_k$ and $L_k=M_k- H_k$.
\begin{lemma}\lab{l:simplexia} 
\bean
\xi&=&O\left(\frac{H_1}{M_1}+\frac{ H_1^3+M_2+L_3}{M_1^2}+\frac{H_1H_2M_2+M_2M_3 }{M_1^3} + \frac{ L_2M_2^2 +   L_2M_2 M_3 }{M_1^4}\right)\\
&&+\, O\left(\frac{M_2^3 L_2+M_2M_3L_4+ L_2L_3H_4}{M_1^5}
\right).
\eean
\end{lemma}
\proof
An upper bound on each $U_k$ is obtained by using either of the two arguments of each min function. Each min function is a function of one or two vertex degrees. If these degrees involved are at least as large as $d_h$, we use the second argument, and otherwise the first. This gives
$$
U_2\le L_2M_2/M_1^2+H_1^2/M_1,
 $$
$$
U_3\le L_4M_3M_2/M_1^4 +H_4L_3L_2/M_1^4+H_3H_1L_3/M_1^3+H_2H_1L_2/M_1^2+H_1^3/M_1,
$$
$$
U_4\le L_3M_2/M_1^2+H_3L_2/M_1^2+H_2H_1/M_1,
$$
$$
U_5\le L_3M_2^2/M_1^4 +H_3L_2^2/M_1^4+H_2H_1L_2/M_1^3+H_1^3/M_1^2,
$$
$$
U_1\le L_3/M_1+H_1.
$$
In $\xi$, we may omit terms that are dominated by others   via inequalities $H_k\le M_k$ and $L_k\le M_k$. These are $ H_1^4/M_1^3\le H_1^3/M_1^2$, and several others involving terms with the same denominators.

The result is
\bean
\xi&=&O\left(\frac{H_1}{M_1}+\frac{ H_1^3+M_2+L_3}{M_1^2}+\frac{H_1H_2M_2+M_2M_3 }{M_1^3}\right)\\
 &&+\,O\left(\frac{H_1^2M_2^2+ M_2^2 L_3+ H_1H_3L_3 +M_2^3+   L_2M_2 H_3 }{M_1^4}+\frac{L_2M_2^3  +M_2M_3L_4+ L_2L_3H_4}{M_1^5}
\right).
\eean
A few more terms can be eliminated as follows.

First, for convenience we replace $L_2M_2H_3/M_1^4$ by $L_2M_2M_3/M_1^4$.  This is tight because $ L_2M_2L_3/M_1^4\le M_2^2L_3/M_1^4$, and the latter is present as a separate error term. 

 Considering the summations in $L_3/L_2$, the ratio of corresponding terms is at most $d_h-2$, whilst in $H_3/H_2$ it is greater. Hence  $L_3/L_2\le d_h-2\le   H_3/H_2$ and consequently $L_3/L_2\le  M_3/  M_2$.
Thus $M_2^2 L_3/M_1^4\le L_2M_2M_3/M_1^4$.

For similar reasons, $ H_1^2M_2^2/M_1^4 \le    H_1H_2M_2/M_1^3$.

By definition $H_1/M_2\le H_1/H_2=O(1/d_h)=O(L_2/L_3)$, and hence
$H_1H_3L_3 \le  L_2M_2H_3    \le L_2M_2M_3$.

Note that $M_2^3=O(H_2^3+L_2M_2^2)$.  Trivially $H_2\le H_1^2\le H_1M_1$, and thus $H_2^3/M_1^4\le H_1H_2^2/M_1^3\le  H_1H_2M_2/M_1^3$, which appears in the bound on   $\xi$.
So we can replace the term $M_2^3/M_1^4$ by $L_2M_2^2/M_1^4$.

The stated bound on $\xi$ follows.\qed

Having $\Delta=\Omega(\sqrt{M_1})$ will permit further simplifications  for an appropriate choice of $h$, as in the following lemma. It is easy to see that, for this value of $h$, these bounds are as tight as that in Lemma~\ref{l:simplexia}.
\begin{lemma}\lab{l:simplexi}
Suppose that $d_h=\Omega(\sqrt{M_1})$ and $d_{h+1}=O(\sqrt{M_1})$ for some  $1\le h\le n-1$. Then
\begin{eqnarray*}
&\mbox{(a)}& \xi = O\left( \frac{H_1 H_2^2  +M_2M_3   }{M_1^3} + \frac{   L_2M_2 M_3 }{M_1^4} + \frac{L_2L_3H_4}{M_1^5}
\right);\\
&\mbox{(b)} & \mbox{if moreover $L_2=\Omega(M_1)$, then the term $M_2M_3/M_1^3$ in (a) can be omitted.}
\end{eqnarray*}
\end{lemma}
\proof
 As $d_h=\Omega(\sqrt{M_1})$, we immediately have $M_2=\Omega(M_1)$ and $M_3=\Omega(M_1)$  and thus $M_2=O(M_3)$. Hence,
$$
\frac{  M_2+L_3}{M_1^2}=O\left(\frac{M_2M_3}{M_1^3}\right), \quad \frac{ L_2M_2^2 }{M_1^4}  =
O\left( \frac{  L_2M_2 M_3 }{M_1^4}
\right).
$$
Moreover, applying Cauchy's inequality shows that $M_2^2=O(M_1M_3)$, so $ L_2 M_2^3/M_1^5 \le L_2M_2 M_3 /M_1^4$.
Hence, the formula in Lemma~\ref{l:simplexia} reduces to
$$
\xi = O\left(\frac{H_1}{M_1}+\frac{ H_1^3 }{M_1^2}+\frac{H_1H_2M_2  +M_2M_3  }{M_1^3} + \frac{   L_2M_2 M_3 }{M_1^4}+ \frac{ M_2M_3L_4+ L_2L_3H_4}{M_1^5}
\right).
$$

Next,    $d_h=\Omega(\sqrt{M_1})$ implies that $H_1^2\ge d_h^2=\Omega(M_1)$ and hence
$H_1/M_1=O(H_1^3/M_1^2)$ which eliminates the first term. It gives moreover that $H_1=O(H_2/\sqrt {M_1})$, and thus
$ H_1^3 /M_1^2 = O(  H_1H_2^2 /M_1^3) = O(H_1H_2M_2 /M_1^3)$, eliminating the second term. Similarly, we get  $H_1H_2L_2/M_1^3 = O(H_2H_3L_2/M_1^4)\le  L_2M_2 M_3/M_1^4$, and thus in the third term $H_1H_2M_2= H_1H_2^2+ H_1H_2L_2$ can be replaced by $H_1H_2^2$.
Finally, $  d_{h+1}=O(\sqrt{M_1})$ implies $L_4\le M_1L_2$, which eliminates $ M_2M_3L_4/M_1^5$.  This gives the bound in part (a). If further we have $L_2=\Omega(M_1)$, then $M_2M_3/M_1^3=O(L_2M_2M_3/M_1^4)$ and part (b) follows.
\qed

 \no {\bf Remark}. It is easy to observe (from the definition of $U_k$ in~\eqn{tildeUs}) that Lemma~\ref{l:simplexi} gives an asymptotically tight bound on $\xi$ if $d_id_j=O(M_1)$ whenever either $i$ or $j$ is at least $h+1$.

\section{Applications}
\lab{sec:applications}

 In this section, we will prove Theorems~\ref{thm:M2}--\ref{thm:heavyPowerLaw} of Section~\ref{sec:results}.
\ss

\noindent{\bf  Power-law density-bounded  degree sequences: proof of Theorem~\ref{thm:powerLaw}}
\smallskip

 Recall  the definition of power-law  density-bounded  degree sequences defined in Section~\ref{sec:results}.  It is easy to see that $\Delta=O(n^{1/\gamma})$ and $M_1=\Theta(n)$. Therefore $\Delta=O( {M_1}^{2/5})$, and so by  Corollary~\ref{cor2:simplexi}(b), the function $\xi$ of Theorem~\ref{thm:general} is $O(M_3^*(M_2^*)^2/M_1^4)$. It is easy to see that $M_k^*=O( n^{(k+1)/\gamma} +n)$ for $k\ge 2$. Hence, $\xi=o(1)$ when $\gamma>5/2$ is fixed. Thus by Theorem~\ref{thm:general} and noting that $1/M_1=O(\sqrt{M_3^*}M_2^*/M_1^2)$,
\be
g({\bf d}) = (1+O(\sqrt{M_3^*}M_2^*/M_1^2))\frac{|\Phi|}{ \prod_{i=1}^n d_i!}\exp\left(-\frac{M_1}{2}+\frac{M_2}{2M_1}+\frac{3}{4}+\sum_{i<j}\log(1+\la_{i,j})\right).\lab{power-lawFormula}
\ee
Next we estimate
$\sum_{i<j}\log(1+\la_{i,j})$.
Taking the Taylor expansion   and   noting that $\sum_{i<j}\sum_{k\ge 4}\la_{i,j}^k=O(M_4^2/M_1^4)$  (since the ratio of the consecutive two terms $\sum_{i<j} \la_{i,j}^{k+1}/\sum_{i<j} \la_{i,j}^{k}$ is $O(M_{k+1}^2/M_k^2 M_1)=O(\Delta^2/M_1)=o(1)$), we have
\bean
\sum_{i<j}\log(1+\la_{i,j})&=&\sum_{i<j} \Big(\la_{i,j}-\frac{\la_{i,j}^2}{2}+\frac{\la_{i,j}^3}{3}\Big)+O(M_4^2/M_1^4).\\
\eean
This we can evaluate using
$$\sum_{i<j}  \la_{i,j}^2=\sum_{i,j} \frac12 \frac{d_i^2d_j^2}{M_1^2}-\sum_i \frac12 \frac{d_i^4}{M_1^2}=(M_2+M_1)^2/2M_1^2 +M_4/4M_1^2  +O(M_3^*/M_1^2)
$$
and so on.  Noting that the error terms $ M_3 ^*/M_1^2+M_2M_3/M_1^3$ are $O(M_3^*(M_2^*)^2/M_1^4)$,
we obtain
\bean
\sum_{i<j}\log(1+\la_{i,j})&=&\frac{M_1}{2}-\frac{M_2}{M_1}-\frac{M_2^2}{4M_1^2}-\frac{3}{4} +\frac{M_3^2}{6M_1^3}+ \frac{M_4}{4M_1^2}-\frac{M_4^2}{8M_1^4}-\frac{M_6}{6M_1^3} +O\left(\frac{M_3^*(M_2^*)^2}{M_1^4}\right).
\eean
Substituting this into~\eqn{power-lawFormula} we obtain the  first formula claimed for $g({\bf d})$.  For the second formula, note that for $\gamma>5/2$ we have $M_4/M_1^2+M_4^2/M_1^4=O(n^{5/\gamma-2})$ and $M_6/M_1^3=O(n^{7/\gamma-3})$, whereas $M_3^*(M_2^*)^2/M_1^4=O(n^{10/\gamma-4})$ if $5/2<\gamma<3$. Hence these terms are all bounded by $O(n^{5/\gamma-2})$ for    $5/2<\gamma<3$.  \qed

\bigskip

\no {\bf Power-law distribution-bounded sequences: proof of Theorem~\ref{thm:iid} }
\ss

Let $Z$ be a random variable with a power-law distribution with parameter $\gamma$, let $p_i=\pr(Z=i)$ and $p_{\ge i}=P(Z\ge i)$.  By the definition of power-law distribution-bounded sequences, the number of vertices with degree at least $i$ is $O(p_{\ge i}n)$.  We may assume that $d_1\ge d_2\ge\cdots\ge d_n$. Then, for every $1\le i\le n$, the number of vertices with degree at least $d_i$ is at least $i$ and so immediately,
$
i=O(p_{\ge d_i} n)=O(n d_i^{1-\gamma})$.
This gives 
\be
d_i=O((n/i)^{1/(\gamma-1)})\quad  \mbox{for}\ 1\le i\le n. \lab{di}
\ee
It follows  that
$M_1=\Theta(n)$ and $M_k=O(n^{k/(\gamma-1)})$ for every fixed integer $k\ge 2$.

Let $\xi$ be defined as in Theorem~\ref{thm:general}. Now we can bound each $U_i$ tightly.  
By~\eqn{di}, for all $i$ and $j$, if $ij\le n^{3-\gamma}$ then $d_id_j=\Omega(n)$ and if $ij>n^{3-\gamma}$ then $d_id_j=O(n)$. This tells how to evaluate the minima in the functions  $U_i$ in~\eqn{tildeUs}.
As an example, we present detailed calculation for a tight bound of $U_4$ only. To simplify the notation, we use $f\preceq g$ to denote $f=O(g)$. 
\bean
U_4&\preceq& \sum_{i=1}^{n^{3-\gamma}} \sum_{j=1}^{n^{3-\gamma}/i} [d_i]_2d_j/n +\sum_{i=1}^{n^{3-\gamma}} \sum_{j=n^{3-\gamma}/i}^{n} [d_i]_3[d_j]_2/n^2 +\sum_{i=n^{3-\gamma}}^n\sum_{j=1}^{n}[d_i]_3[d_j]_2/n^2\\
&\preceq& \sum_{i=1}^{n^{3-\gamma}} \sum_{j=1}^{n^{3-\gamma}/i} \frac{(n/i)^{2/(\gamma-1)} (n/j)^{1/(\gamma-1)}}{n}+\sum_{i=1}^{n^{3-\gamma}} \sum_{j=n^{3-\gamma}/i}^{n} \frac{(n/i)^{3/(\gamma-1)} (n/j)^{2/(\gamma-1)}}{n^2}\\
&&+\frac{M_2}{n^2}\sum_{i=n^{3-\gamma}}^n (n/i)^{3/(\gamma-1)}.
\eean
From here, it is straightforward to obtain
\[
\frac{U_4M_2}{M_1^2}=O(n^{\frac{(3-\gamma)(\gamma-2)+8-3\gamma}{\gamma-1}})=O(n^{\frac{2+2\gamma-\gamma^2}{\gamma-1}}),
\]
which is $o(1)$ when $\gamma>1+\sqrt{3}$.
We neglect the calculations of the other terms in $\xi$ as the approach is similar and these terms are dominated by $U_4M_2/M_1^2$ when $\gamma>1+\sqrt{3}$. This completes the proof of the theorem. \qed
 
\bigskip

\no {\bf Using $M_1$ and $M_2$ alone: proof of Theorem~\ref{thm:M2}}\ss

We can assume $M_2\ge 1$ since otherwise $d_i=1$ for all $i$ and the theorem is is true with zero error term.
Then we have $M_2\ge [\Delta]_2$ and $\Delta=O(M_2^{1/2})$. Choose $h$ to be the minimum integer for which   $d_{h+1}\le \sqrt{M_1}$. If $h\ge 1$, we can easily bound $L_2$ and $H_2$ by $M_2$; $M_3$ and $H_3$ by $O(M_2^{3/2})$ (since $H_3\le M_3\le \Delta M_2$); $L_3$ by $O(\sqrt{M_1} M_2)$ and $L_4$ by $O( M_1  M_2)$; and $H_4$ by $M_2^2$ (since $H_4\le \Delta^2 H_2=O(M_2^{2}))$.
 Note that  $M_2\ge H_2\ge (d_h-1) H_1 = \Omega(M_1^{1/2}H_1)$ by the choice of $h$,   and so $H_1 = O(M_2 /M_1^{1/2})$.  Define $\xi$ as in Theorem~\ref{thm:general}.  Then, by Lemma~\ref{l:simplexi}(a) (and by noting that $M_2=\Omega(\sqrt{M_1})$ as $h\ge 1$),
 $\xi=O(M_2^4/M_1^{9/2})$.  If $h=0$ (i.e.\ $\Delta\le \sqrt{M_1}$), we can bound $M_3$ by $M_2^{3/2}$ and then Corollary~\ref{cor2:simplexi}(b) gives
$\xi=O((M_2^{3/2}+M_1)(M_2^2+M_1^2)/M_1^4) = O( M_2^4/M_1^{9/2}+M_2^{3/2}/M_1^2+1/M_1)$ (as the terms apart from $1/M_1$ are dominated by $M_2^4/M_1^{9/2}$  when $M_2=\Omega(M_1)$ and by $M_2^{3/2}/M_1^2 $  otherwise).
  Thus, $\xi=o(1)$ as long as $M_2=o(M_1^{9/8})$, and Theorem~\ref{thm:M2} follows (with its redefinition of $\xi$) by Theorem~\ref{thm:general}.
\qed

\bigskip

\noindent{\bf Bi-valued sequences: proof of Theorem~\ref{thm:bidegrees}}
\smallskip

Now $M_i=[\Delta]_i\ell+[\delta]_i(n-\ell)$ for every integer $i\ge 1$. It is easy to see that $M_1=\Theta(\Delta\ell+\delta n)$. Since $\delta\ge 3$, we have $M_1=O(M_2)$ and $M_1=O(M_3)$.
Define $\xi$ as in Theorem~\ref{thm:general}. We first show that $\xi=o(1)$ as long as either (a) or (b) holds.

By the hypotheses in (a), $\Delta=O(\sqrt{M_1})$. Applying Corollary~\ref{cor3:simplexi},
\[
\xi=\Theta\left(\frac{(\Delta^3\ell+\delta^3n)(\Delta^4\ell^2+\delta^{4}n^2)}{\Delta^4\ell^4+\delta^4n^4}\right)=O\left(\frac{\Delta^7\ell^3+\Delta^3\delta^4\ell n^2+\delta^7n^3}{\Delta^4\ell^4+\delta^4n^4}\right)=o(1),
\]
as it is easy to verify that $\Delta^4\delta^3\ell^2 n=O(\Delta^7\ell^3+\delta^7n^3)$.
This proves part (a).

Now we prove part (b). By our assumptions, $\Delta^5\ell^3=o(\delta^3 n^3)$ and $\Delta=\Omega((\delta n)^{1/2})$. We first show that
 $\Delta \ell< \delta n$. Suppose not, then
 \[
\frac{\Delta^5\ell^3}{\delta^3n^3}\ge \Delta^2=\Omega(\delta n),
\]
contradicting the assumption that $\Delta^5\ell^3/\delta^3n^3=o(1)$.
Thus $\Delta\ell<\delta n$ and immediately $M_1=\Theta(\delta n)$. So $\Delta=\Omega(\sqrt{\delta n})$ implies that $\Delta=\Omega(\sqrt{M_1})$.
  Applying Lemma~\ref{l:simplexi} with $h=\ell$, we have (using $\delta\le \Delta\le n$)
\bean
\xi&=&O\left(\frac{\Delta^5\ell^3}{\delta^3 n^3}+\frac{(\Delta^3\ell+\delta^3 n)\delta^2 n (\Delta^2\ell+\delta^2 n)}{\delta^4n^4}+\frac{\Delta^4\ell \delta^5 n^2}{\delta^5n^5}\right)\\
&=&O\left(\frac{\Delta^5\ell^3}{\delta^3n^3}+\frac{\Delta^5\ell^2}{\delta^2n^3}+\frac{\delta^3}{n}+\frac{\Delta^3\ell }{n^2}\right)=o(1).
\eean
 We have now shown that under any condition of (a,b), $\xi=o(1)$. It is easy to see that both $M_3/M_1^2$ and $1/M_1$ are dominated by $\xi$ in each case. So the theorem follows by Theorem~\ref{thm:general}.\qed
\smallskip

\no {\bf Remark}. We have obtained as strong a result as if we had evaluated the expression for $\xi$  in Theorem~\ref{thm:general} directly rather than using the results of Section~\ref{sec:boundxi}. This follows by the remark after Lemma~\ref{l:simplexi}, and by noting  that in (b) we can assume $\delta\Delta<M_1$ (since $M_1=\Theta(\delta n)$ and $ \Delta<n$). The results in (a) and (b) are similarly tight.

\bigskip

\no {\bf Long-tailed power-law degree sequences: proof of Theorem~\ref{thm:heavyPowerLaw}}\ss

  Choose $h$ to be the minimum integer for which $d_{h+1}<n^\alpha$. If $h=0$, the degrees are uniformly bounded, which is a case treated in~\cite{BC}. However, the error term there is only $o(1)$. Instead, we are done  by~\cite[Theorem 4.6]{M2}, where the error term is $O(1/n)$ which is clearly $O(\sqrt{\xi})$, with $\xi$ as defined in the theorem statement.    Otherwise, $d_h=\Omega(n^{\alpha})=\Omega(\sqrt{M_1})$,  since $\alpha>1/2$.  Moreover,  by part (a) of the definition of these degree sequences, any component that is less than $n^{\alpha}$ is bounded. So $d_{h+1}= O(\sqrt{M_1})$, and we can apply Lemma~\ref{l:simplexi}.  
   It is easy to verify that $\Delta=O(n^{\alpha+\beta/\gamma})$ and, in the notation of Lemmas~\ref{l:simplexi} and~\ref{l:simplexia},
$H_1= O(n^{\alpha+\beta} )$ for $\gamma> 2$,  $H_1= O(n^{\alpha+\beta}\log n )$ for $\gamma= 2$,  and $H_1= O(n^{\alpha+2\beta/\gamma})$ for $1<\gamma<2$. By our assumption on $\beta$, it is easy to verify that $H_1=o(n)$ always. For every fixed $k\ge 2$, $H_k=O(n^{k\alpha+(k+1)\beta/\gamma})$ since $\gamma<3$, and $L_k=O(n)$, and moreover $L_1=\Theta(n)$.  This implies that $M_1=\Theta(n)$. Now define $\xi$ as in Theorem~\ref{thm:general}.
 By Lemma~\ref{l:simplexi} and using $\alpha>1/2$, it is easy to check that
  \[
  \xi=\left\{
  \begin{array}{ll}
  O(n^{5\alpha+\beta+6\beta/\gamma-3} ) & \mbox{if $2< \gamma<3$}\\
    O(n^{5\alpha+\beta+3\beta-3}\log n )  & \mbox{if $\gamma=2$}\\
  O(n^{5\alpha+8\beta/\gamma-3}) & \mbox{if $1<\gamma<2$}.
  \end{array}
  \right.
  \]
  By the assumption on $\beta$, we have $\xi=o(1)$.  As $\alpha>1/2$ by our assumption, the bound on $\xi$ presented above obviously dominates $1/n$. The theorem now follows by Theorem~\ref{thm:general}. \qed

\ss

\no
{\bf Acknowledgement} 

We are grateful to Remco van der Hofstad  for a communication on power law distribution-bounded sequences.


\begin{thebibliography}{99}

\bibitem{ACL} W. Aiello, F. Chung and L. Lu,
Random evolution in massive graphs, {\em 42nd IEEE Symposium on Foundations of Computer Science} (2001), {510-�519}.

\bibitem{ACL2} W. Aiello, F. Chung, and L. Lu,
A random graph model for power law graphs, {\em Experiment. Math.} {\bf 10} (2001), no. 1, {53-�66}.

\bibitem{AB} R. Albert and A. Barab\'{a}si, Statistical mechanics of complex networks, {\em Rev. Mod. Phys}, (2002).



\bibitem{BC} E.A. Bender and E.R. Canfield, The asymptotic number of
labeled graphs with given degree sequences, {\em J. Combinatorial
Theory Ser.\ A} {\bf 24} (1978),  {296--307}.


\bibitem{BP} M. Bogu\~{n}\'{a}, F. Papadopoulos and D. Kriouko, Sustaining the Internet with hyperbolic mapping, {\em Nature
communications} {\bf 1} (62), 2010.


\bibitem{B} B. Bollob\'{a}s, A probabilistic proof of an asymptotic formula
for the number of  labelled regular graphs, {\em Europ. J.
Combinatorics} {\bf 1} (1980), {311--316}.


\bibitem{BRST} B. Bollob\'{a}s, O. Riordan, J. Spencer and G. Tusn\'{a}dy,
The degree sequence of a scale-free random graph process,
{\em Random Structures Algorithms} {\bf 18} (2001), no. 3, {279-�290}.


\bibitem{BCDR} C. Borgs, J. Chayes, C. Daskalakis and S. Roch,
First to market is not everything: an analysis of preferential attachment with fitness {\em STOC'07 -- Proceedings of the 39th Annual ACM Symposium on Theory of Computing}, (2007), {135-�144}.

\bibitem{CF} C. Cooper and A. Frieze,
A general model of web graphs,
{\em Random Structures Algorithms} {\bf 22} (2003), no. 3, {311�-335}.

\bibitem{CL} F. Chung and L. Lu,
The volume of the giant component of a random graph with given expected degrees,
{\em SIAM J. Discrete Math.} {\bf 20} (2006), no. 2, {395�-411}.


\bibitem{ER} P. Erd\H{o}s and A. R\'{e}nyi,
On random graphs. I.
{\em Publ. Math. Debrecen} 6 1959, {290--297}.


\bibitem{FFF}  	M. Faloutsos, P. Faloutsos and 	
	C. Faloutsos, On power-law relationships of the Internet topology, {\em ACM SIGCOMM Computer Communication Review}
{\bf 29} 4 (1999), {251 -- 262}.

\bibitem{GPP}  	L. Gugelmann, K. Panagiotou and U. Peter, Random hyperbolic graphs: degree sequence and clustering, {\em Proceeding
ICALP'12 Proceedings of the 39th international colloquium conference on Automata, Languages, and Programming} {\bf II},
{573--585}.

\bibitem{J} S. Janson,
The probability that a random multigraph is simple,
{\em Combin. Probab. Comput.} {\bf 18}, no. 1--2, (2009) {205--225}.

\bibitem{J2} S. Janson, The probability that a random multigraph is simple, II. arXiv:1307.6344


\bibitem{K} R. Kleinberg, Geographic routing using hyperbolic space, {\em
Proceedings of the 26th IEEE International Conference on Computer Communications, INFOCOM�07}
1902�-1909.

\bibitem{M2}   B.D. McKay, Asymptotics for symmetric 0-1
matrices with prescribed row sums, {\it Ars Combinatoria} {\bf 19A} (1985),
15--25.


\bibitem{MW2} B.D. McKay and N.C. Wormald, Asymptotic enumeration by degree
sequence of graphs with degrees $o(\sqrt{n})$, {\em Combinatorica} {\bf
11} (1991), {369--382}.

\bibitem{M} M. Mitzenmacher,
A brief history of generative models for power law and lognormal distributions,
{\em Internet Math} {\bf 1} (2004), no. 2, {226-�251}.

\bibitem{PR} C.H. Papadimitriou and D. Ratajczak, On a conjecture related to geometric routing,
{\em Theoretical Computer
Science} {\bf 344} 1, (2005), {3-�14}.


\bibitem{PKBV} F. Papadopoulos, D. Krioukov, M. Bogu\~{n}\'{a} and A. Vahdat, Greedy forwarding in dynamic scale-free networks
embedded in hyperbolic metric spaces, {\em Proceedings of the 29th IEEE International Conference on Computer
Communications, INFOCOM�10},
{2973-�2981}.
\bibitem{regsurvey} N.C. Wormald, Models of random regular graphs, {\it
Surveys in Combinatorics, 1999}, London Mathematical Society Lecture Note
Series {\bf 267} (J.D. Lamb and D.A. Preece, eds) Cambridge University
Press,  Cambridge,
pp.\ {239--298}, 1999.

\end{thebibliography}
\end{document}